\documentclass[12pt, amssymb]{article}

\def\sect#1{\advance\sectionnum by 1\section{#1}}
\newcount\sectionnum \sectionnum=-1
\def\thesection{\the\sectionnum}

\newtheorem{theorem}{Theorem}[section]
\newtheorem{lemma}[theorem]{Lemma}

\newtheorem{definition}[theorem]{Definition}
\newtheorem{proposition}[theorem]{Proposition}
\newtheorem{corollary}[theorem]{Corollary}
\newtheorem{remark}[theorem]{Remark}

\usepackage{amsmath, amssymb, latexsym}

\begin{document}

\title{\bf  Topologies on the group of Borel automorphisms of a
standard Borel space}
\author{{\bf S.~Bezuglyi}\thanks{Supported in part by  CRDF grant
UM1-2546-KH-03}\\ Institute for Low Temperature Physics, Kharkov, Ukraine\\
\and{\bf A.H.~Dooley}\\ University of New South Wales, Sydney, Australia\\
\and {\bf J.~Kwiatkowski}\thanks{Supported in part by KBN grant 1P 03A 038
26}\\ College of Economics and Computer Science, Olsztyn, Poland}

\date{}

\maketitle

\begin{abstract}
The paper is devoted to the study of topologies on the group $Aut(X,{\cal
B})$ of all Borel automorphisms of a standard Borel space $(X, {\cal B})$.
Several topologies are introduced and all possible relations between them
are found. One of these topologies, $\tau$, is a direct analogue of the
uniform topology widely used in ergodic theory. We consider the most
natural subsets of $Aut(X,{\cal B})$ and find their closures. In
particular, we describe closures of subsets formed by odometers, periodic,
aperiodic, incompressible, and smooth automorphisms with respect to the
defined topologies. It is proved that the set of periodic Borel
automorphisms is dense in $Aut(X,{\cal B})$ (Rokhlin lemma) with respect to
$\tau$. It is shown that the $\tau$-closure of odometers (and of rank 1
Borel automorphisms) coincides with the set of all aperiodic automorphisms.
For every aperiodic automorphism $T\in Aut(X,{\cal B})$, the concept of a
Borel-Bratteli diagram is defined and studied. It is proved that every
aperiodic Borel automorphism $T$ is isomorphic to the Vershik
transformation acting on the space of infinite paths of an ordered
Borel-Bratteli diagram. Several applications of this result are given.

\end{abstract}

\newcommand{\La}{\Lambda}
\newcommand{\Om}{\Omega}
\newcommand{\om}{\omega}
\newcommand{\e}{\varepsilon}
\newcommand{\al}{\alpha}
\newcommand{\vp}{\varphi}
\newcommand{\G}{\Gamma}
\newcommand{\g}{\gamma}
\newcommand{\N}{{\mathbb N}}
\newcommand{\T}{\theta}
\newcommand{\De}{\Delta}
\newcommand{\de}{\delta}
\newcommand{\s}{\sigma}
\newcommand{\E}{{\cal E}}
\newcommand{\A}{{\cal A}}
\newcommand{\B}{{\cal B}}
\newcommand{\M}{{\cal M}}
\newcommand{\bs}{(X,{\cal B})}
\newcommand{\Aut}{Aut(X,{\cal B})}
\newcommand{\h}{Homeo(\Om)}
\newcommand{\pos}{{\mathbb R}^*_+}
\newcommand{\R}{{\mathbb R}}
\newcommand{\la}{\lambda}
\newcommand{\Z}{{\mathbb Z}}
\newcommand{\ap}{{\cal A}p}
\newcommand{\per}{{\cal P}er}
\newcommand{\inc}{{\cal I}nc}

\newpage

\sect{Introduction}

The study of topologies on the group of transformations of an underlying
space has a  long history. Some of the early results in the area are the
classical results of J.~Oxtoby and S.~Ulam on the typical dynamical
behavior of homeomorphisms which preserve a measure [OU, O]. Traditionally,
this circle of problems has attracted attention in various areas of
dynamical systems, notably, in measurable and topological dynamics, where
it is important for many applications to understand what kind of
transformations is typical for certain dynamics. Of course, this problem
assumes that a topology is defined on the group of all transformations. The
best known results concerning ergodic, mixing, and weakly mixing
automorphisms of a measure space were obtained by P.~Halmos and
V.A.~Rokhlin (see, e.g. [H], [R]). Many results on approximation of
automorphisms of a measure space can be found in the book by I.~Cornfeld,
S.~Fomin, and Ya.~Sinai [CFS]. The recent book by S.~Alpern and V.S.~Prasad
[AP] develops the Oxtoby-Ulam approach and contains new results on
approximation of homeomorphisms of compact and non-compact manifolds.

Motivated by ideas used in measurable dynamics, we first started studying
topologies in the context of Cantor minimal systems in [BK1] and [BK2],
which we refer to as Cantor dynamics, but very soon it became clear that
our approach could also be used for Borel automorphisms of a standard
Borel space (Borel dynamics). The main goal of our two papers (see also
[BDK]) is to study the global properties of topologies on the group of all
Borel automorphisms $\Aut$ of a standard Borel space $\bs$ and the group
of all homeomorphisms $\h$ of a Cantor set $\Om$. Although we define and
study several topologies on $\Aut$, only two of them, the uniform $(\tau)$
and weak $(p)$ topologies, are considered as basic since they are
analogous to the topologies which are well-known in ergodic theory. Let us
recall their definitions for automorphisms of a measure space.

Let $(X, \B,\mu)$ be a standard measure space. On the group
$Aut(X,\B,\mu)$ of all non-singular automorphisms of $X$, the uniform and
weak topologies are defined by metrics $d_u(S,T) =\mu(\{x\in X\ :\ Sx\neq
Tx\}) $ and $ d_w(S,T) = \sum_n 2^{-n}\mu(SA_n\ \De\ TA_n)$, respectively,
where $S,T\in Aut(X,\B,\mu)$ and $(A_n)$ is a countable collection of
Borel sets which is dense in $\B$. These concepts have turned out to be of
crucial importance in ergodic theory. As far as we know, the first deep
results on the topological properties of $Aut(X,\B,\mu)$ with respect to
$d_u$ and $d_w$ appeared in the pioneering papers by P.~Halmos [H1] and
V.A.~Rokhlin [R] where the concept of approximation of automorphisms was
introduced in abstract ergodic theory. Later, this concept was developed
in many papers where the notion of approximation was considered in various
areas of ergodic theory. Probably the most famous application of these
ideas in measurable dynamics is the Rokhlin lemma, a statement on
approximation of an aperiodic automorphism by periodic ones. However,
there have been many other applications, too numerous to list here.

It seems rather curious but, to the best our knowledge, these topologies
have not so far been systematically studied in the context either of
topological or Borel dynamics. However, we should mention the interesting
paper [GlW] by E.~Glasner and B.~Weiss where the Rokhlin property is
considered  for homeomorphisms of a compact metric spaces with respect to
the topology of uniform convergence. For a standard Borel space, we do not
have a fixed Borel measure on the underlying space (in contrast to
measurable dynamics). Therefore, if we want to extend the definitions of
topologies generated by metrics $d_u$ and $d_w$ to the group $\Aut$, then
we have to take into account the set ${\cal M}_1(X)$ of {\it all} Borel
probability measures on $\bs$. Roughly speaking, we say that two Borel
automorphisms $S$ and $T$ from $\Aut$ are close in the uniform topology
$\tau$ on $\Aut$ if for {\it any} measure $\mu \in {\cal M}_1(X)$ the set
where $S$ and $T$ are different is small in the measure $\mu$. To define
the topology $p$, which is treated as an analogue of $d_w$, we observe
that if the symmetric difference of two Borel sets is arbitrarily small
with respect to {\it any} $\mu\in {\cal M}_1(X)$, then these sets must
coincide. Thus, a $p$-neighborhood of $T$ is formed by those $S$ for which
$SB = TB$ where $B$ is a given Borel set (see Section 1 for strict
definitions).

In [BK1, BK2], we first gave a definition of the uniform topology $\tau$
analogous to $d_u$. Our main interest in those papers was focused on full
groups and normalizers of Cantor minimal systems. In particular, we showed
that the full group of a minimal homeomorphism is closed in $\tau$. The
paper [BDK] is devoted to the study of topologies on the group $\h$ of all
homeomorphisms of a Cantor set $\Om$. We have tried to make the content of
the current paper parallel to that of [BDK] by considering similar
topologies and similar problems. The comparison of all results from these
papers would take up too much place. We  mention only that on a Cantor set
$\Om$ one can study both dynamics, Borel and Cantor, and answer some
questions about their interplay because $\h$ is obviously a subset of
$Aut(\Om,\B)$. In particular, we can consider in both cases the topology
of uniform convergence generated by the metric $D(S,T) = \sup_{x\in \Om}
d(Sx,Tx) + \sup_{x\in \Om} d(S^{-1}x,T^{-1}x)$ on the groups $\h$ and
$Aut(\Om,\B)$. Notice that $D$ coincides with the topology $p$ defined on
$\h$ by clopen sets only [BDK]. This topology is extremely useful in the
study of homeomorphisms of a Cantor set. In Borel dynamics, the topologies
$p$ and $D$ are obviously inequivalent. Nevertheless, it seems to be
interesting to study topological properties of the metric space
$(Aut(\Om,\B), D)$ keeping in mind the parallel theory for Cantor
dynamics.

The outline of the paper is as follows. We consider more systematically
the definitions and properties of various natural topologies on $\Aut$.
Actually, we study simultaneously a collection of Hausdorff topologies:
$\tau$ (the analogue of uniform topology), $\tau'$ (which is equivalent to
$\tau$), $\tau''$ (which is weaker than $\tau$), $p$ (which is the direct
analogy of the weak topology in ergodic theory and which is mostly useful
in the context of Cantor dynamics), $\tilde p$ (which is equivalent to
$p$), and $\overline{p}$ (which is weaker than $p$). We consider also the
topologies $\tau_0$ (which is weaker than $\tau$) and $p_0$ (which is
equivalent to $p$) as natural modifications of $\tau$ and $p$. They all
(except $\overline{p}$) make $\Aut$ into a topological group.

The initial part of the paper is devoted to discovering all possible
relations between these topologies (see Theorem \ref{T1.2} and its proof
in  Section 3). In Section 1, we discuss various topological properties of
the group $\Aut$ and its subsets. For example, we describe convergent
sequences of automorphisms from $\Aut$ and show that $(\Aut, p)$ is a
zero-dimensional topological space (Corollary \ref{C1.7}). The group
$\Aut$ has a normal subgroup $Ctbl(X)$ which consists of automorphisms
with at most countable support. It turns out that $Ctbl(X)$ is closed with
respect to the above topologies. This fact allows us to study topologies
on the quotient group $\widehat{Aut}(X,\B)= \Aut/Ctbl(X)$. This kind of
identification of Borel automorphisms is analogous to that usually used in
ergodic theory. In Sections 2 and 4, we study the following classes of
Borel automorphisms: periodic, aperiodic, smooth, incompressible, and of
rank 1 (the latter includes odometers). It is proved that periodic
automorphisms are dense in $\Aut$ with respect the topology $\tau$. This
allows us to prove a version of the Rokhlin lemma (Theorem \ref{T2.3}).
Remark that the problem of periodic approximation of an aperiodic Borel
automorphism has been also studied by M.~Nadkarni [N] and B.~Weiss [W]. As
an immediate consequence of these results, we obtain that the set $\ap$ of
aperiodic automorphisms is nowhere dense in $(\Aut, \tau)$.  We also prove
that the full group $[T]$ of a Borel automorphism $T\in \Aut$ is closed
with respect to the all topologies. We consider the set of incompressible
automorphisms, ${\cal I}nc$,  consisting of those aperiodic automorphisms
which admit an invariant Borel probability measure and prove that ${\cal
I}nc$ is a closed nowhere dense subset of $\ap$ with respect to $p$. It is
shown that the $\tau$-closure of rank 1 Borel automorphisms coincides with
that of odometers (Theorem \ref{T2.8}).  In the last section, we introduce
the concept of Borel-Bratteli diagrams in the context of Borel dynamics.
To do this, we use a remarkable result on existence of a  vanishing
sequence of markers [BeKe, N]. It is well known that Bratteli diagrams
play a very important role in the study of minimal homeomorphisms of a
Cantor set [HPS, GPS]. Similarly to Cantor dynamics, we show that every
aperiodic Borel automorphism is isomorphic to the Vershik transformation
acting on the space of infinite paths of a Borel-Bratteli diagram. Several
applications of this result are given. In particular, we prove that the
set of odometers is $\tau$-dense in aperiodic Borel automorphisms. We also
discuss properties of such diagrams for automorphisms of compact and
locally compact spaces. We believe that this approach to the study of
Borel automorphisms will lead to further interesting developments in this
area.

To conclude, we would like to make two remarks. Firstly,  Borel
automorphisms of a standard Borel space (or, more general, countable Borel
equivalence relations) have been extensively studied in many papers during
last decade. We refer to numerous works of A.~Kechris, G.~Hjorth,
M.~Foreman, M.~Nadkarni, B.~Weiss and others. More comprehensive
references can be found, for instance, in [BeKe, Ke2, Hjo, FKeLW, N].
Secondly, it is impossible to discuss in one paper all interesting
problems related to topologies on $\Aut$. We consider the current paper as
the first step in the study of topological  properties of $\Aut$. Our
primary goal is to display a wealth of new topological methods in Borel
dynamics.
\\

Throughout the paper, we use the following standard {\bf notation}:

\begin{itemize}

\item $(X, \B)$ is a standard Borel space with the $\sigma$-algebra of
Borel sets $\B= \B(X)$; $\B_0$ is the subset of $\B$ consisting of
uncountable Borel sets.

\item $Aut(X, \B)$ is the group of all one-to-one Borel automorphisms of $X$
with the identity map ${\mathbb I} \in Aut(X, \B)$.

\item $\ap$ is the set of all aperiodic Borel automorphisms.

\item $\per$ is the set of all periodic Borel automorphisms.

\item ${\cal M}_1(X)$ is the set of all Borel probability measures on
$(X,\B)$. Let ${\cal M}^c_1(X)$ denote the subset of ${\cal M}_1(X)$ formed
by continuous (non-atomic) Borel probability measures.

\item $\de_x$ is the Dirac measure at $x \in X$.

\item $E(S,T)=\{ x\in X \ |\ Tx\ne Sx\} \cup \{x\in X \ |\ T^{-1}x\ne
S^{-1}x\}$ where $S,T \in \Aut$.

\item $B(X)$ is the the set of Borel real-valued bounded functions, $B(X)_1
= \{f\in B(X)\ |\ \Vert f \Vert := \sup\{|f(x)| : x\in X\} \le 1\}$, $\mu
(f)=\int_X f\,d\mu$ where $f\in B(X)$, $\mu\in {\cal M}_1(X)$.

\item $\mu\circ S(A):=\mu (SA) $ and $\mu \circ S(f): =\int_{X}f\,d(\mu
\circ S) = \int_{X}f(S^{-1}x)\,d\mu (x)$ where $S\in \Aut$.

\item $A^c := X \setminus A$ where $A\in \B$.

\item the term \textquotedblleft automorphism"  means a Borel automorphism
of $\bs$; we deal with Borel subsets of $X$ only.

\item $\Om$ is a Cantor set.

\item $\h$ is the group of all homeomorphisms of $\Om$.

\end{itemize}

\sect{Topologies on $\Aut$}

\setcounter{equation}{0}

{\bf 1.1. Definition of topologies on ${\Aut}$}. In this section, we define
topologies on $\Aut$ and establish their main properties.

\begin{definition}\label{D1.1}  The  topologies $\tau,\tau',
\tau'',\ p,\ \tilde p$, and $\overline{p}$ on $\Aut$ are defined by the
bases of neighborhoods ${\cal U},\ {\cal U}',\ {\cal U}'',\ {\cal W},\
\widetilde{\cal{W}}$, and $\overline{\cal W}$, respectively. They are:
${\cal U} = \{U(T; \mu_1,...,\mu_n; \e)\}$, ${\cal U'} = \{U'(T;
\mu_1,...,\mu_n; \e)\},\ {\cal U''} = \{U''(T; \mu_1,...,\mu_n; \e)\}$,
${\cal W} = \{W(T; F_1,...,F_k)\}$, $\widetilde{\cal{W}}= \widetilde{W}(T;
f_1,...,f_n;\e)$, and $\overline{\cal W} =\{\overline W(T; F_1,...,F_k;\
\mu_1,...,\mu_n;\ \e)\}$ where
\begin{equation}\label{1.1}
U(T; \mu_1,...,\mu_n; \e)= \{ S\in \Aut\ |\ \mu_i(E(S,T))<\e ,\
i=1,...,n\},
\end{equation}
\begin{equation}\label{1.2}
U'(T; \mu_1,...,\mu_n; \e )= \{ S\in \Aut \ |\ \sup_{F\in
{\cal B}}\mu_i(TF\ \De \ SF)<\e,\ i=1,...,n\},
\end{equation}
\begin{equation}\label{1.3}
\begin{array}{ll}
U''(T; \mu_1,...,\mu_n; \e )& = \{ S\in \Aut\ |\ \sup_{f\in B(X)_1} \vert
\mu_i\circ S(f)-\mu_i\circ T(f)\vert\\
& <\e, \ i=1,...,n\},
\end{array}
\end{equation}
\begin{equation}\label{1.4}
W(T; F_1,...,F_k) = \{S\in \Aut\ |\ SF_i = TF_i, \ i= 1,...,k\},
\end{equation}
\begin{equation}\label{1.4a}
\widetilde{W}(T;f_1,...,f_m;\e) = \{S\in \Aut\ |\ \Vert f_i \circ T^{-1} -
f_i \circ S^{-1}\Vert <\e, i=1,...,m\}
\end{equation}
\begin{equation}\label{1.5}
\begin{array}{ll}
\overline{W}(T; (F_i)_1^k; (\mu_j)_1^n; \e)&= \{S\in \Aut\ |\ \mu_j(SF_i\
\De TF_i) + \mu_j(S^{-1}F_i \De T^{-1}F_i)\\
 &< \e,\ i=1,...,k; j=1,...,n\}.
\end{array}
\end{equation}
In all the above definitions $ T\in \Aut $, $ \mu_1,...,\mu_n \in {\cal
M}_1(X),\ F_1,...,F_k \in \B$, $f_1,...,f_m\in B(X)$, and $\e>0 $.
\end{definition}

If the set $E_0(S,T) = \{x\in X : Sx\neq Tx\}$ were used in (\ref{1.1})
instead of $E(S,T)$, then we would obtain the topology equivalent to
$\tau$.

It is natural to define also two further topologies, which are similar to
$\tau$ and $p$, by considering only continuous measures and uncountable
Borel sets.

\begin{definition}\label{D1.1a} The topologies $\tau_0$ and $p_0$ on $\Aut$
are defined by the bases of neighborhoods ${\cal U}_0 =
\{U_0(T;\nu_1,...,\nu_n;\e)\}$ and ${\cal W}_0 = \{W_0(T;A_1,...,A_n)\}$,
respectively, where $U_0(T;\nu_1,...,\nu_n;\e)$ and $W_0(T;A_1,...,A_n)$ are
defined as in (\ref{1.1}) and (\ref{1.4}) with $\nu_i \in {\cal M}_1^c(X)$
and $A_i\in \B_0,\ i=1,...,n$.
\end{definition}

Obviously, $\tau_0$ is not stronger than $\tau$ and $p_0$ is not stronger
than $p$.

Note that the topology $\tau$ was first introduced in [BK1] where,
motivated by ergodic theory, we called it the {\it uniform topology}. We
defined $p$ in [BK2] in the context of homeomorphisms of Cantor sets. In
this section, we use a number of results about these topologies which are
proved lately. Namely, the following theorem is proved in Section 3.

\begin{theorem}\label{T1.2}  $(1)$ The topologies $\tau$ and $\tau'$
are equivalent. \\ $(2)$ The topology $\tau\ (\sim \tau')$ is strictly
stronger than $\tau''$. \\ $(3)$ The topology $\tau$ is strictly stronger
than $\overline{p}$. \\ $(4)$ The topology $\tau$ is strictly stronger than
$\tau_0$.\\ $(5)$ The topology $p$ is equivalent to $\tilde p$.\\ $(6)$ The
topology $p\ (\sim\tilde p)$ is equivalent to $p_0$.\\ $(7)$ The topology
$p$ is strictly stronger than $\overline{p}$. \\ $(8)$ The topology $\tau$
is not comparable with $p$ and the topology $\tau''$ is not comparable with
$\overline{p}$ and $\tau_0$.
\end{theorem}

We will also introduce in Section 3 two auxiliary topologies $\tilde \tau$
and $\overline{\tau}$ equivalent to $\tau''$ which will allow us to have a
more convenient description of $\tau''$ (see Definition \ref{D3.3} and
Proposition \ref{P3.4}). In particular, the topology $\overline{\tau}$ is
defined by neighborhoods
\begin{equation}\label{tau''}
\overline V(T; \mu_1,...,\mu_n; \e ) =  \{S\in \Aut\ \vert\ \sup_{F\in
{\cal B}}\vert\mu_j(TF) - \mu_j(SF)\vert <\e,\ j=1,...,n\}
\end{equation}
where $\mu_1,...,\mu_n\in {\cal M}_1(X)$.

Given an automorphism $T$ of $(X,{\cal B})$, we can associate a linear
unitary operator $L_T$ on the Banach space $B(X)$ by $(L_Tf)(x) =
f(T^{-1}x)$. Then the topology $\tilde p$ is induced on $\Aut$ by the
strong operator topology on bounded linear operators of $B(X)$. Theorem
\ref{T1.2} asserts equivalence of $p$ and $\tilde p$. It is well known
that the weak topology in ergodic theory can be also defined in a similar
way using the strong topology on linear bounded operators of a Hilbert
space. This observation is a justification of the name {\it
\textquotedblleft weak topology"} which will be used to refer to $p$.

\begin{proposition}\label{P1.3}  ${\cal U}, {\cal U'}, {\cal U''}, {\cal
U}_0, {\cal W}, \widetilde{\cal W}, {\cal W}_0$, and $\overline{\cal W}$
are bases of Hausdorff topologies $\tau, \tau', \tau'', \tau_0,\ p,\ \tilde
p,\ p_0$, and $\overline{p}$, respectively. $\Aut$ is a topological group
with respect to $\tau,\ \tau', \tau'', \tau_0,\ p,\ \tilde p,\ p_0$, and
$\Aut$ is not a topological group with respect to $\overline{p}$.
\end{proposition}

\noindent {\it Proof}. The first statement is clear for $\tau,\ \tau',
\tau_0,\ p,\ \tilde p,\ p_0$ and can be immediately deduced from the
definition of neighborhoods (\ref{1.1}), (\ref{1.2}), (\ref{1.4}),
(\ref{1.5}). We need to check it for $\tau''$ and $\overline{p}$ only. Let
$S$ be a Borel automorphism taken in a given neighborhood $U''(T;\
\mu_1,...,\mu_n;\ \e)$. We will show that there exists $U''_0(S;\
\nu_1,...,\nu_n;\ \de)\subset U''(T;\ \mu_1,...,\mu_n;\ \e)$. To see this,
take $\nu_i=\mu_i, \ \de =\e -c$, where
$$
c=\max_{1\le i\le n}\{ \sup_{f\in B(X)_1} \vert \mu_i\circ T(f)-
\mu_i\circ S(f)\vert \}.
$$
If $R\in U''_0$, then we get
$$
\begin{array}{lll}
& \sup_{f\in B(X)_1} \vert \mu_i\circ T(f)-\mu_i\circ R(f)\vert\\
\\
\le &
\sup_{f\in B(X)_1}\vert \mu_i\circ T(f)-\mu_i\circ S(f)\vert
+\sup_{f\in B(X)_1}\vert \mu_i\circ S(f)-\mu_i\circ R(f)\vert\\
\\
 < & c+\e -c=\e,
\end{array}
$$
i.e. $R \in U''(T;\ \mu_1,...,\mu_n;\ \e )$.

To see that $\{\overline{\cal W}\}$ is a base of neighborhoods, we may use a
slight modification of the above argument. For $S\in \overline{W}(T;
F_1,...,F_k; \mu_1,...,\mu_n; \e)$, take
$$
c=\max_{i,j}[\mu_j(SF_i\ \De\ TF_i) + \mu_j(S^{-1}F_i\ \De\ T^{-1}F_i)].
$$
It is easily seen from (\ref{1.5}) that
$$
\overline{W}(S; F_1,...,F_k; \mu_1,...,\mu_n; \de)\subset \overline{W}(T;
F_1,...,F_k; \mu_1,...,\mu_n; \e)
$$
where $\de = \e - c$.

Now we prove that $\Aut$ is a topological group with respect to $\tau\ \sim\
\tau',\ \tau_0,\ p\ \sim p_0\ \sim\ \tilde p$, and $\tau''$.

Consider first $(\Aut, \tau)$. In the case of the topology $\tau_0$ the
proof is similar. Let $S,T\in \Aut$. We need the following facts:

(i) $ U(T; \mu_1,...,\mu_n;\e) = U(T^{-1}; \mu_1,...,\mu_n;\e)^{-1}$;

(ii) $\{x\in \Om : Sx\neq Tx\} \subset \{x\in \Om : Sx\neq Rx\}\cup
\{x\in \Om : Rx\neq Tx\}$.

\noindent
Indeed, (i) follows from the relation $E(S,T) = E(S^{-1},T^{-1})$ and
(ii) is checked straightforward.

By (i), the map $T\mapsto T^{-1}$ is continuous. To prove that $ (S,T)
\mapsto ST $ is also continuous, we show that for any neighborhood $ U_{ST}
= U(ST; \mu_1,...,\mu_n;\e) $ there exist neighborhoods $U_{S} = U(S;
\nu_1,...,\nu_k;\e_1) $ and $ U_{T} = U(T; \sigma_1,...,\sigma_m;\e_2) $
such that $ U_SU_T\subset U_{ST} $. Take $\e_1 =\e_2 =\e/4 $, $k=m= 2n$,
$(\nu_1,...,\nu_{2n}) = (\mu_1,...,\mu_n, \ \mu_1\circ T^{-1},...,\mu_n\circ
T^{-1}) $, and $ (\sigma_1,...,\sigma_{2n}) =(\mu_1,...,\mu_n, \mu_1\circ
S,...,\mu_n\circ S) $. Let $P\in U_S$, $ Q\in U_T $. It follows from (ii)
that
$$
\begin{array}{lllll}   &E(PQ, ST)  =   \{ x : PQx \neq STx\} \cup  \{ x
:( PQ)^{-1}x \neq (ST)^{-1}x\}\\
\\
  \subset &\{ x : PQx \neq PTx\} \cup \{ x : PTx \neq STx\}\\
& \cup  \{ x : Q^{-1}P^{-1}x \neq Q^{-1}S^{-1}x\} \cup \{ x :
Q^{-1}S^{-1}x \neq T^{-1}S^{-1}x\}\\
\\
=&  \{ x : Qx \neq Tx\} \cup  T^{-1}(\{ x : Px \neq Sx\})\\
 &\cup  \{ x :P^{-1}x \neq S^{-1}x\} \cup S(\{ x : Q^{-1}x \neq
T^{-1}x\}).
\end{array}
$$
Then, for any $i=1,...,n$,
\begin{eqnarray*}
\mu_i(E(PQ,ST)) &\le &\mu_i (\{ x : Qx \neq Tx\})+\mu_i\circ T^{-1}(\{ x
: Px \neq Sx\})\\
\\
&+ &\mu_i( \{ x :P^{-1}x \neq S^{-1}x\})+\mu_i\circ S(\{ x :
Q^{-1}x \neq T^{-1}x\}) <  \e.
\end{eqnarray*}
Thus, $PQ\in U_{ST}$.

The proof of continuity of $(T,S) \mapsto TS$ and $T\mapsto T^{-1}$ in
the space $(\Aut, p)$ is straightforward.

To complete the proof of the theorem, we will show that $\Aut$ is a
topological group with respect to $\tau''$. Again choose $S,T \in \Aut$ and
let $U''(ST) = U''(ST, \mu_1,...,\mu_n;\e)$ be a given neighborhood.
Consider $U''(S) = U''(S;\mu_1,...,\mu_n; \e/2)$ and $U''(T) =
U''(T;\mu_1\circ S,...,\mu_n\circ S;\e/2)$. Then for $R\in U''(S), Q\in
U''(T)$ we have
$$
\sup_{f\in B(X)_1}\vert \mu_i\circ RQ(f) - \mu_i\circ ST(f)\vert
$$
$$
\le
\sup_{f\in B(X)_1}\vert \mu_i\circ RQ(f) - \mu_i\circ SQ(f)\vert +
\sup_{f\in B(X)_1}\vert \mu_i\circ SQ(f) - \mu_i\circ ST(f)\vert
$$
$$
 = \sup_{f'=L_Qf\in B(X)_1}\vert \mu_i\circ R(f') -
\mu_i\circ S(f')\vert + \sup_{f\in B(X)_1}\vert (\mu_i\circ S)Q(f) -
(\mu_i\circ S)T(f)\vert < \e.
$$

We remark that $(\Aut, \overline{p})$ is not a topological group as one can
prove that the product $(T,S) \mapsto TS$ is not continuous in the topology
$\overline{p}$. To see this, we show that there exists a
$\overline{p}$-neighborhood $\overline{W}_0$ of the identity such that in
any neighborhoods $\overline{W}_1(\mathbb I; (E_i); (\nu_j);\e)$ and
$\overline{W}_2(\mathbb I; (D_i); (\la_j);\de)$ of the identity one can
find automorphisms $R$ and $Q$,respectively, such that $RQ\notin
\overline{W}_0$. Take $\overline{W}_0= \overline{W}_0(\mathbb I; F;
\de_{x_0}; 1/2)$. Note that if $Sx_0 \notin F$, then $S\notin
\overline{W}_0$. Find an automorphism $Q\in \overline{W}_2$ such that the
point $Qx_0$ is of measure zero with respect to all $\nu_j$. Finally, find
some $R\in \overline{W}_1$ such that $RQx_0\notin F$. \hfill{$\square$}

\begin{remark}\label{rem} {\rm  (1) Let $(X, {\cal B})$ and $(Y,{\cal C})$ be
two standard Borel spaces (hence, they are Borel isomorphic).  It can
easily be seen that $(\Aut, top(X))$ and $(Aut(Y,{\cal C}), top(Y))$ are
homeomorphic where $top$ is any of the topologies from Definitions
\ref{D1.1} and \ref{D1.1a}.

(2) Let $(X,d)$ be a compact metric space. In this case, we can consider
the group $Aut(X,\B)$ of Borel automorphisms and its subgroup $Homeo(X)$ of
homeomorphisms of $X$. Define for $S,T\in \Aut$ the topology of uniform
convergence generated by metric\footnote{In fact, the metric $D$ can be
considered on $\Aut$ if $X$ is a totally bounded metric space. Many of
results concerning the metric $D$ which are proved below for a compact
metric space can be generalized to totally bounded spaces.}
\begin{equation}\label{DD}
D(S,T) = \sup_{x\in X} d(Sx,Tx) + \sup_{x\in X} d(S^{-1}x,T^{-1}x).
\end{equation}
Then $(\Aut, D)$ is a complete metric space and $Homeo(X)$ is closed in
$\Aut$.  When $X= \Om$ is a Cantor set, we studied thoroughly the Polish
group $(Homeo(\Om),D)$ in [BDK]. It is not hard to see that in Cantor
dynamics the topology on $\h$ generated by $D$ is equivalent to the
topology $p$ defined by clopen sets only.  We note that, in contrast to
(1), the topology generated by $D$ on $\Aut$ depends, in general, on the
topological space $(X,d)$. Nevertheless, we think it is worth to study the
topological properties of $\Aut$ and $Homeo(X)$ for a fixed compact (or
Cantor) metric space  $X$ because one can compare in this case these
properties for the both groups.}
\end{remark}

Let $Ctbl(X)$ be defined as the subset of $\Aut$ consisting of all
automorphisms with countable support, that is $T\in Ctbl(X)
\Leftrightarrow (E(S, \mathbb I)\ {\rm is \ at\ most\ countable}$. One can
show that $Ctbl(X)$ is a normal subgroup, closed with respect to the
topologies which we have defined (see Lemma \ref{ctbl} below). Therefore
$\widehat{Aut}(X,{\cal B})= \Aut/Ctbl(X)$ is a Hausdorff topological group
with respect to the quotient topology. Note that the quotient group
$\widehat{Aut}(X,{\cal B})$ was first considered in [Sh] where the
simplicity of this group was proved. Considering elements from
$\widehat{Aut}(X,{\cal B})$, we identify Borel automorphisms which are
different on a countable set. The class of automorphisms equivalent to a
Borel automorphism $T$ we will denote by the same symbol $T$ and write
$T\in \widehat{Aut}(X,{\cal B})$. This corresponds to the situation in
ergodic theory when two automorphisms are also identified if they are
different on a set of measure 0. We studied topological properties of the
group $\widehat{Aut}(X,{\cal B})$ in [BM].

\begin{lemma}\label{ctbl} The normal subgroup $Ctbl(X)\subset \Aut$ is
closed with respect to the all topologies from Definitions \ref{D1.1} and
\ref{D1.1a}.
\end{lemma}

\noindent{\it Proof}. By Theorem \ref{T1.2}, it is sufficient to prove the
statement of the lemma for the topologies $\tau'', \tau_0$ and
$\overline{p}$. Notice that we gave in [BM] a direct proof of the fact
that $Ctbl(X)$ is closed in $\tau$ and $p$. The case of the topology
$\tau_0$ is similar to that of $\tau$.

Suppose that $S \in \overline{Ctbl(X)}^{\tau''} \setminus Ctbl(X)$. Then
there exists an uncountable Borel set $F$ such that $SF\cap F =\emptyset$.
Let $\mu$ be a continuous Borel probability measure such that $\mu(F) = 1$
and let $f_0(x) = \chi_F(x)$. The $\tau''$-neighborhood $U''(S) =
U''(S;\mu;1/2)$ consists of automorphisms $T$ satisfying the condition
$$
\sup_{f\in B_1(X)} |\mu\circ S(f) - \mu\circ T(f)| < 1/2.
$$
By assumption, there exists $T_0 \in Ctbl(X)\cap U''(S)$. Then $|\mu\circ
S(f_0) - \mu\circ T_0(f_0)| < 1/2$. On the other hand, it can be easily
seen that $\mu\circ T_0(f_0) =1$ and  $\mu\circ S(f_0) = 0$. This leads to
a contradiction.

The fact that $Ctbl(X)$ is closed in $\overline{p}$ is proved in the same
way. We leave the details to the reader. \hfill$\square$
\\

Notice also that $\overline{Ctbl(X)}^D = Ctbl(X)$ when $(X, d)$ is a
compact metric space. Indeed, if $(T_n)_{n\in \N} \subset Ctbl(X)$, then
there exists a countable set $C\subset X$ such that $T_nx = x \ \forall
x\in X\setminus C$. Let $D(T_n,T) \to 0 \ (n\to \infty)$. We obtain  $Tx =
x$ for every  $x\in X\setminus C$, i.e. $T\in Ctbl(X)$.

Let $\pi$ be the natural projection from $\Aut$ to $\widehat{Aut}(X,{\cal
B})$. Lemma \ref{ctbl} allows us to define the quotient topologies
$\hat\tau = \pi(\tau)$, $\hat\tau_0 = \pi(\tau_0)$, and $\hat p = \pi(p)$
on $\widehat{Aut}(X,{\cal B})$. It turns out that $\hat\tau$-neighborhoods
are defined by continuous  measures and $\hat p$-neighborhoods are defined
by uncountable Borel sets. In particular, this means that $\hat\tau$ is
equivalent to $\hat\tau_0$. In [BM], the following proposition was proved.

\begin{proposition}\label{coun} Given a $\hat\tau$-neighborhood  $\hat U=
\hat U(T;\mu_1,...,\mu_n;\e)$ and a $\hat p$-neighborhood $\hat W= \hat
W(T;F_1,...,F_m)$, there exist neighborhoods $\hat
U_0(T;\nu_1,...,\nu_n;\e) = U_0(T;\nu_1,...,\nu_n;\e)Ctbl(X)$ and $\hat
W_0(T;B_1,...,B_m) = W_0(T;B_1,...,B_m)Ctbl(X)$ in $\hat\tau$ and $\hat p$,
respectively, such that $\hat U_0 \subset \hat U,\ \hat W'\subset \hat W$
where $\nu_1,...,\nu_n\in {\cal M}_1^c(X)$ and $B_1,...,B_m \in \B_0$.
\end{proposition}

\noindent {\bf 1.2. Properties of the topologies}. We will now discuss
some properties of the  topologies which we introduced above. In
particular, we consider convergent sequences with respect to these
topologies.

\begin{remark}\label{R1.4} {\rm (1) We recall some facts about the
uniform topology $\tau$ from [BK1].

\noindent
(a) $T_n\stackrel{\tau}{\longrightarrow}S \iff \forall x \in X\
\exists n(x)\in \N$
such that $\forall n >n(x), \ T_nx = Sx$;

\noindent
(b)
\begin{eqnarray*}
T_n\stackrel{\tau}{\longrightarrow}S & \iff & \forall \mu \in {\cal
M}_1(X), \ \mu(E(T_n,S))\to 0\\
& \iff & \forall x\in X, \ \de _x(E(T_n,S)) \to 0;
\end{eqnarray*}

\noindent (c) $(\Aut,\tau)$ is a complete nonseparable topological group in
the sense that any Cauchy sequence of Borel automorphisms converges to a
Borel automorphism;

\noindent
(d) For a Cantor set $\Om$, $Homeo(\Om )$ is not closed in $Aut(\Om,
\B(\Om))$ with respect to $\tau$. (In fact, it is shown in [BDK] that
$Homeo(\Om )$ is dense in  $Aut(\Om, \B(\Om))$ in $\tau$).

(2) Note that
\begin{equation}\label{tau0}
T_n\stackrel{\tau_0}{\longrightarrow}T \iff [\forall \mu\in \M_1^c(X),\
\mu(E(T_n,T)) \to 0,\ n\to \infty].
\end{equation}
It is not hard to see that (\ref{tau0}) holds if the set
\begin{equation}\label{cnt}
C = \bigcap_{n\in \N}\bigcup_{m\geq n} E(T_m,T)
\end{equation}
is countable. Indeed, if $C_n = \bigcup_{m\geq n}E(T_m,T)$, then $\mu(C_n)
\to 0$ \ $\forall \mu\in \M_1^c(X)$. Hence $\mu(E(T_m,T)) \to 0$ as $m\to
\infty$.

Observe that condition (\ref{cnt}) is not necessary for
$\tau_0$-convergence. This means that even a weaker form of (1)(a) does not
hold. To see this, let us consider the following example. Set $X = [0,1)$
and denote by $(\xi_n)$ the refining sequence of partitions of $X$ into the
intervals $A_n(i) = [i2^{-n}, (i+1)2^{-n}),\ i =0,1,...,2^n -1,\ n\in \N$.
Let $\mu$ be a continuous Borel probability measure on $X$. Then $\forall
\e >0$\ $\exists N_\e\in \N$ such that $\mu(A_n(i)) < \e$\ $\forall n >
N_\e$.

Now let $T_{n}(i)$ be a Borel automorphism of $X$ such that $E(T_{n}(i),
\mathbb I) = A_{n}(i)$. Then $\mu(E(T_{n}(i), \mathbb I) \to 0$ as $n\to
\infty$. By (\ref{tau0}), this sequence of automorphisms converges to the
identity map in $\tau_0$. On the other hand, for any $x\in [0,1)$ the
property $T_{n}(i)x \neq x$ holds for infinitely many automorphisms from
the sequence $(T_n(i))$.

(3) $(\Aut, p)$ is a complete nonseparable topological group. Note that it
follows from Theorem \ref{T1.2} and (1) that $(\Aut, \tau'')$ and $(\Aut,
\overline{p})$ are also complete spaces.

(4) Let $T$ be a Borel automorphism and let $U(T;\mu_1,...,\mu_n;\e)$ be a
$\tau$-neighborhood of $T$. Consider $\nu = n^{-1}(\mu_1 +\dots + \mu_n)$.
Then it can be easily shown that $U(T;\ \nu;\ n^{-1}\e) \subset
U(T;\mu_1,...,\mu_n;\e)$. This means that we can work with a single measure
instead of a finite collection of measures when it is more convenient.}
\end{remark}

\begin{remark}\label{R1.5} {\rm We notice the following  three simple
properties of the topology $p$.\\ (1) Without loss of generality, we may
assume that $p$ is generated by neighborhoods $W(T;F_1,...,F_n)$ where
$(F_1,...,F_n)$ is a partition of $X$. \\ (2) Let $W(T;F_1,...,F_n)$ be
given, then for every $S\in W(T;F_1,...,F_n)$ one has
$W(S;F_1,...,F_n)=W(T;F_1,...,F_n)$. It follows from this observation that
$W(\mathbb I; F_1,...,F_n)$ is an open subgroup of $(\Aut, p)$.\\ (3) For
any $T\in \Aut$, we have $W(T^{-1}; TF_1,...,TF_n)^{-1}= W(T;F_1,...,F_n)$.}
\end{remark}

\begin{proposition}\label{P1.6}  The sets $W(T;F_1,...,F_n)$ are closed
in $\Aut$ with respect to the topologies $\tau,\tau',\tau'', p$ and
$\overline{p}$.
\end{proposition}

\noindent {\it Proof}. We begin with two simple observations. Firstly, the
proposition will be proved if we can show that every set $W(T;F)$ is
closed. Secondly, if the proposition holds for $\tau''$ and
$\overline{p}$, then it holds for the other topologies because they all
are stronger than either $\tau''$ or $\overline{p}$.
\smallskip

(1) Let us  first consider $\overline{p}$. We show that $W(T;F)^c$ is open
in $\overline{p}$. Suppose $S\in W(T;F)^c$, that is $SF \neq TF$. Then we
have two cases:

(i) $E:= TF \setminus SF \neq \emptyset$,

(ii) $TF \setminus SF =\emptyset,$ that is $TF \subset SF$.\\ In case (i),
define a neighborhood $\overline{W}(S; F;\mu;1/2)$ with a measure $\mu$
concentrated on $E$, i.e. $\mu(E) = 1, \mu (E^c)=0$ (one can take $\mu =
\de_x$ with some $x\in E$). Let $R\in \overline{W}(S;F;\mu;1/2)$, then
$\mu(RF) = \mu (RF\ \De\ SF) <1/2$ since $\mu(SF) = 0$. On the other hand,
$\mu (TF\setminus RF) \ge \mu(TF) - \mu(RF) > 1/2$ since $TF\supset E$. This
means that $TF \neq RF$ and $R \in W(T;F)^c$.

In case (ii), we take $\mu$ such that $\mu(SF \setminus RF) =1$. Consider
again an automorphism $R\in \overline{W}(S;F;\mu;1/2)$. We have
$$
1/2 > \mu(SF\ \De\ RF) = \mu (RF\setminus SF)+ \mu(SF\setminus RF) =
\mu((SF\setminus TF)\setminus RF)
$$
and therefore
$$
\mu((SF\setminus TF)\cap RF)= \mu(SF\setminus TF) - \mu((SF\setminus
TF)\setminus RF) > 1/2.
$$
Thus, $(SF\setminus TF)\cap RF \neq \emptyset$ and $RF \neq TF$, i.e. $R\in
W(T;F)^c$.
\smallskip

(2) Show that $W(T;F)^c$ is open in $\tau''$. In fact, we will consider the
topology $\overline{\tau}$ equivalent to $\tau''$ (see (\ref{tau''}) and
Proposition \ref{P3.4}). Let $S\in W(T;F)^c$, then $SF\neq TF$ and we have
the above cases (i) and (ii). In case (i), take a
$\overline{\tau}$-neighborhood $\overline{V}(S;\mu;1/2)$ as in
(\ref{tau''}) where $\mu$ is concentrated on $E$. Then $\mu(SF) = 0$ and
for $R\in \overline{V}(S;\mu;1/2)$ we have that
$$
\mu(RF) = \vert\mu(RF) -\mu(SF)\vert < 1/2.
$$
Therefore, $\mu(TF\setminus RF) > 1/2$ and $R\in W(T;F)^c$.

In case (ii), we choose $\mu$ supported on $SF\setminus TF$. Then, by the
method of (1), we see that for $R\in \overline{V}(S;\mu;1/2)$ the set
$(SF\setminus TF)\cap RF$ is non-empty since its measure is greater than
1/2. Therefore, $RF\neq TF$ and $R\in W(T;F)^c$. \hfill$\square$
\\

Observe that Proposition \ref{P1.6} does not hold for $\tau_0$. One can
easily show that $W(T;\{x\})$ is not closed with respect to $\tau_0$.

\begin{corollary}\label{C1.7}  $(\Aut, p)$ is a $0$-dimensional
topological space.
\end{corollary}

\noindent {\it Proof}. This follows from Proposition \ref{P1.6}.
\hfill$\square$
\\

Now we consider convergent sequences in each of our topologies. Of course,
the topologies are not defined by convergent sequences but it is useful for
many applications to know criteria of convergence.

\begin{remark}\label{R1.8} {\rm (1) If a sequence $(T_n)$ of Borel automorphisms
converges to $S$ in $\tau''$, then for any measure $\mu\in {\cal M}_1(X)$,
$$
\vert \mu\circ S(f)-\mu\circ T_n(f)\vert \to 0
$$
uniformly in $f\in B(X)_1$ as $n\to\infty$. Since $\tau''$ is equivalent to
$\overline{\tau}$ (Proposition 3.4), the above condition is equivalent to
the following one: $T_n\stackrel{\tau''}{\longrightarrow}S$ if and only if
$$
\vert \mu (SF)-\mu (T_nF)\vert \to 0
$$
uniformly in $F\in \B$.

(2) $(T_n)$ converges to $S$ in $p$ if and only if for any Borel set $F$,
$T_nF = TF$ for sufficiently large $n = n(F)$. In particular, $F$ can be a
point from $X$. Therefore, we see that $p$-convergence implies
$\tau$-convergence by Remark \ref{R1.4}.

(3) It follows directly from (\ref{1.5}) that $T_n\stackrel{\overline p}
{\longrightarrow} S$ if and only if $\forall \mu \in {\cal M}_1(X) \ \forall
F \in \B$}
\begin{equation}\label{1.6}
 \mu(T_nF\ \De\ SF)+\mu (T^{-1}_nF\ \De\ S^{-1}F)\to 0.
\end{equation}
\end{remark}

In fact, one can prove the following criterion of
$\overline{p}$-convergence.

\begin{proposition}\label{P1.9}  $T_n\stackrel{\overline p}
{\longrightarrow} S$ if and only if $\forall F \in \B$,
\begin{equation}\label{1.7}
SF=\limsup_{n\to\infty} T_nF,\ \ \ S^{-1}F= \limsup_{n\to\infty} T_n^{-1}F,
\end{equation}
where
$$
\limsup_{n\to\infty} F_n = \bigcup_m \bigcap_{n>m}F_n.
$$
\end{proposition}

\noindent {\it Proof}. We assume for simplicity that $S= {\mathbb I}$. The
general case is proved similarly. To prove (\ref{1.7}) we remark that for
any $x\in X$ and $F\in \B$, the convergence $T_n\stackrel{\overline p}
{\longrightarrow} {\mathbb I}$ implies that
$$
\de_x (T_nF\ \De \ F)+\de_x (T^{-1}F\ \De\ F) \to 0
$$
as $n\to \infty$. This means that if $x\in F$, then there exists
$n_0=n_0(x,F)$ such that $x\in T_nF$ and $x\in T^{-1}_nF$ for all $n>n_0$.
We have proved that $F \subset \bigcup_m \bigcap_{n>m} T_nF,\ F\subset \bigcup_m
\bigcap_{n>m} T^{-1}_nF$. In fact, these inclusions are equalities. Indeed, if
we assume that there exists $x_0\in F^c=X\setminus F$ with $x_0\in
\bigcap_{n>m}T_nF$ for some $m$, then we have a contradiction to the fact that
$x_0$ also belongs to $\bigcup_k \bigcap_{n>k}T_nF^c$. Thus, (\ref{1.7}) holds.

Conversely, let $E_m=\bigcap_{n>m}T_nF$ and $\bigcup_m E_m=F$. Since
$E_m\subset E_{m+1}$, we know that for any measure $\mu \in {\cal M}_1(X),\
\mu E_m\to \mu F\ (m\to \infty)$. Remark that $E_m\subset T_nF$ for all
$n>m$. Therefore $E_m=E_m\cap T_nF\subset F\cap T_nF\subset F$. Thus we
obtain $\mu (F\cap T_nF)\to \mu F,\ n\to \infty$. Similarly $\mu (F\cap
T^{-1}_nF)\to \mu F$. By (\ref{1.6}), the proof is complete.
\hfill{$\square$}

\begin{proposition}\label{P1.10} Suppose $(T_n)$ is a sequence
of Borel automorphisms and let $S\in \Aut$. Then
$$
(T_n \stackrel{p}{\longrightarrow} S)\ \Rightarrow \
(T_n\stackrel{\tau}{\longrightarrow} S)\ \iff\
(T_n\stackrel{\tau''}{\longrightarrow}S) \ \iff \ (T_n\stackrel{\overline
p} {\longrightarrow}S).
$$
\end{proposition}

\noindent {\it Proof}. We will consider, for simplicity, the case $S =
{\mathbb I}$. We note that if $T_n \stackrel{p}{\longrightarrow} {\mathbb
I}$, then for every $x \in X$ the sequence $(T_n)$ eventually gets into
$W({\mathbb I}; \{x\})$, that is $T_nx = x$ for sufficiently large $n$. It
follows that $T_n\stackrel{\tau}{\longrightarrow}{\mathbb I}$ (see Remark
\ref{R1.4}).

By Theorem \ref{T1.2}, $\tau$ is strictly stronger than $\tau''$. To prove
the second implication, we need to verify only that if a sequence of Borel
automorphisms $(T_n)$ converges in $\tau''$, then it also converges in
$\tau$.

To see this, we assume that $T_n \stackrel{\tau''}\longrightarrow {\mathbb
I}$ when $n\to \infty$. Then for any measure $\mu\in {\cal M}_1(X)$, we have
$$
\sup_{f\in B(X)_1} \vert \mu(f)-\mu\circ T_n(f)\vert \to 0,\qquad n\to
\infty.
$$
Thus,
\begin{equation}\label{1.8}
\left|\int_X(f(T_n^{-1}x) - f(x))d\mu\right| \to 0
\end{equation}
uniformly in $f$ as $n\to\infty$. Take $\mu=\de_{x_0}$ in (\ref{1.8}). Then
$\forall \e>0\ \exists N=N(\e,x_0)$ such that $\forall n > N$ and $\forall
f\in B(X)_1$, one has
\begin{equation}\label{1.9}
\vert f(T_n^{-1}x_0) - f(x_0)\vert < \e.
\end{equation}
To prove that $T_n\stackrel{\tau}\longrightarrow {\mathbb I}$, it suffices
to verify that $T_nx_0=x_0$ for $n$ sufficiently large (see Remark \ref{R1.4}).
To obtain a contradiction, we assume that for any $N$ there exist a point
$x_0$ and $n_0 > N$ such that $T_{n_0}x_0 \neq x_0$. Take a Borel set $F$
of $x_0$ such that $T_{n_0}x_0 \notin F$. Then there exists a function
$f_0\in B(X)_1$ such that $f_0(x_0) =1$ and $f_0(T_{n_0}x_0) = 0$. This
contradicts (\ref{1.9}).

As above we will show that $\overline{p}$-convergence of $(T_n)$ implies
$\tau$-convergence. If $T_n\stackrel{\overline p} {\longrightarrow}{\mathbb
I}$, then for every $x\in X$, $T_n\in \overline{W}({\mathbb I}; \{x\};
\de_x; 1/2)$ when $n$ is sufficiently large. In other words, $T_nx = x$. It
proves that $T_n\stackrel{\tau} {\longrightarrow}{\mathbb
I}$.\hfill$\square$


\sect{Approximation by periodic and aperiodic automorphisms}

\setcounter{equation}{0}

\noindent {\bf 2.1. Periodic approximation of Borel automorphisms}. Here we
focus on the study of periodic and  aperiodic automorphisms. We will show
that for every Borel aperiodic automorphism $T$ of $\bs$ there exists a
sequence of periodic Borel automorphisms that converges to $T$ in the
uniform topology $\tau$. In fact, this result was proved by Nadkarni in [N]
although he did not consider topologies on $\Aut$. We reproduce the main
part of Nadkarni's proof here because it will be used below. We will also
find the closures of some natural classes of automorphisms.

Recall some standard definitions. For $T\in \Aut$, a point $x\in X$ is
called {\it periodic} if there exists $n\in \N$ such that $T^nx=x$. The
smallest such $n=n(x)$ is called the {\it period} of $T$ at $x$. Given $T$,
the space $X$ can be partitioned into a disjoint union of Borel
$T$-invariant sets $X_1, X_2,...,X_\infty$ where $X_n$ is the set of points
with period $n$, and $X_\infty$ is the set where $T$ is aperiodic. Such a
partition related to an automorphism $T$ will be called canonical. If
$X_\infty =\emptyset$, then $T$ is called {\it pointwise periodic}, $T \in
\per$. Denote by $\per_n(x)$ the set of all automorphisms which have period
$n$ at $x$. By definition, $T\in \per_n$, the set of all Borel
automorphisms of period $n$, if $X_n = X$. In other words,
\begin{equation}\label{per_n}
\per_n =\bigcap_{x\in X} \per_n(x).
\end{equation}
We say that $T\in \per_0$ if there exists $N\in \N$ such that $P^Nx=x,\
x\in X$. This means that $X$ is a finite union of some sets
$X_{n_1},...,X_{n_k}$. Obviously, $\per_0$ is a proper subset of $\per$.
Finally, if $X = X_\infty$, then $T$ is called {\it aperiodic}, $T\in \ap$.

\begin{proposition}\label{P2.1} $(1)$ For any
$n\in \N$, the set $\per_n(x)\ (x\in X)$ is clopen with respect to all
topologies from Definition \ref{D1.1}.\\
$(2)$ \ $\overline{\per_0}^\tau = \overline{\per}^\tau$.\\
$(3)$ \ $\ap$ and $\per_n\ (n\in \N)$ are closed with respect to the all
topologies.
\end{proposition}
{\it Proof}. (1) As mentioned above, to show that $\per_n(x)$ is closed for
the all topologies, it suffices to do this for $\tau''$ and $\overline{p}$.

By Proposition \ref{P3.4}, $\tau''$ is equivalent to the topology
$\overline{\tau}$ whose neighborhoods are defined by $\overline V(T;\
\mu_1,...,\mu_n;\ \e ) =  \{S\in \Aut\ \vert\ \sup_{F\in {\cal
B}}\vert\mu_j(TF) - \mu_j(SF)\vert <\e,\ j=1,...,n\}$ (see Definition
\ref{D3.3}). Let $R\in \overline{\per_n(x)}^{\overline{\tau}}$. Let
$\overline{V}(R) = \overline{V}(R; \de_{Rx},...,\de_{R^nx}; 1/2)$. Then,
there exists $P\in \overline{V}(R)\cap \per_n(x)$. This means that taking
one-point sets $\{x\},...,\{R^{n-1}x\}$, we obtain that
$$
|\de_{R^ix}(R(R^{i-1}x)) - \de_{R^ix}(P(R^{i-1}x))| < 1/2,\ \ i =1,...,n.
$$
It follows that $Rx = Px$,\dots ,$R^nx = P^nx = x$, i.e. $R\in \per_n(x)$.

Let us show that $\per_n(x)$ is closed in $\overline{p}$. Take $R\in
\overline{\per_n(x)}^{\overline{p}}$ and consider the
$\overline{p}$-neighborhood $\overline{W}(R) = \overline{W}(R;
\{x\},...,\{R^{n-1}x\}; \de_{Rx},...,\de_{R^nx}; 1/2)$. There is an
automorphism $P\in \overline{W}(R) \cap \per_n(x)$. The inequalities
$\de_{Rx}(Rx \De Px) < 1/2,\dots, \de_{R^nx}(R(R^{n-1}x) \De P(R^{n-1}x)) <
1/2$ imply that $Rx = Px,\dots,R^nx = P^nx =x$.

To check that $\per_n(x)$ is open with respect to the all topologies, it
suffices again to show this for $\tau''$ and $\overline{p}$ only (Theorem
\ref{T1.2}) . This fact follows from the following observation:  if $P\in
\per_n(x)$, then $\overline{V}(P; \de_{Px},...,\de_{P^{n}x}; 1/2)$ and
$\overline{W}(P; \{x\},...,\{P^{n-1}x\}; \de_{Px},...,\de_{P^{n}x};1/2)$
are subsets of $\per_n(x)$.
\medskip

(2) Let now $T\in \per$. We can construct a sequence $(P_n)\subset \per_0$
converging to $T$ in $\tau$. As above, consider the partition $(X_i:
i=1,2,...)$ of $X$ corresponding $T$, i.e. $T$ has period $i$ on $X_i$. Let
$E_n = \bigcup_{i\le n}X_i$. Then $E_n \subset E_{n+1}$ and $X = \bigcup_n
E_n$. Define $P_nx = Tx$ if $x\in E_n$ and $P_nx = x$ if $x \in X\setminus
E_n$. Clearly, $P_n \stackrel{\tau}\longrightarrow T$.
\medskip

(3) The set $\per_n$ is closed by (1) and (\ref{per_n}). It is clear that
$$
\ap = \Aut \setminus \left(\bigcup_{n\in \N}\bigcup_{x\in X}
\per_n(x)\right).
$$
Hence $\ap$ is closed.

\hfill$\square$
\\

Let $Orb_T(x)$ denote the $T$-orbit of $x\in X$.  Recall  the
definition of the {\it full group} $ [T] $ generated by $T\in \Aut$:
$$
[T] = \{ \g \in \Aut\ \vert\ \g x \in Orb_T(x),\ \forall x\in X\}.
$$
Then every $\g \in [T]$ defines a Borel function $m_\g : X \to \Z$
such that $\g x = T^{m_\g (x)}x,\ x\in X $. Thus, every $\g \in
[T]$ defines a countable partition of $X$ into Borel sets $A_n =
\{ x\in X: m_\g(x) =n\},\ n\in \N$.

It is obvious that if $T\in \widehat{Aut}(X,{\cal B})$, then one can also
define the full group $[T]$ as a subgroup of $\widehat{Aut}(X,{\cal B})$.

\begin{proposition}\label{[T]} $(1)$ The full group $[T]\ (T\in \Aut)$
is closed in $\Aut$ with respect to the topologies from Definition
\ref{D1.1}.\\
$(2)$ The full group $[T]\ (T\in \widehat{Aut}(X,{\cal B}))$ is closed in
$\widehat{Aut}(X,{\cal B})$ with respect to $\hat\tau$.\footnote{We do not
consider here other quotient topologies on $\widehat{Aut}(X,{\cal B})$.}
\end{proposition}
{\it Proof}. (1) It is not hard to prove this result directly for the
topologies $\tau$ and $p$. To prove the proposition for all our
topologies, it is sufficient to check that  $[T]$ is closed in $\tau''
\sim \overline{\tau}$ and $\overline{p}$.

Assume that there exists $S\in \overline{[T]}^{\overline{\tau}} \setminus
[T]$. Then one can find a point $y\in X$ such that $Sy \neq T^ny$ for all
$n\in \Z$. Let $\mu$ be an atomic probability measure supported by $\{T^ny
: n\in \Z\}$ such that $\mu(\{y\}) =1/2$. Then the
$\overline{\tau}$-neighborhood $\overline{V}(S) = \overline{V}(S;\mu;1/4)$
contains an automorphism $\g\in [T]$. Hence for any Borel set $F$ we have
that $|\mu(SF) - \mu(\g F)| < 1/4$. For $F = \{\g^{-1}y\}$, we have a
contradiction.

The proof for $\overline{p}$ is similar. We observe only that the
$\overline{p}$-neighborhood $\overline{W}(S;\{y\}; \de_{Sy}; 1/2)$ cannot
meet $[T]$ where $S$ and $y$ as above.

(2) Assume that $R\in \overline{[T]}^{\hat\tau}\setminus [T]$ where $T\in
\widehat{Aut}(X,{\cal B})$. Then any $\hat\tau$-neighborhood $\hat U(R)$
contains an element $\g$ from $[T]$. Since $R$ is not in $[T]$, the Borel
set $A = \bigcap_{n\in \Z}E(R,T^n)$ is uncountable. Let $\mu$ be a
continuous measure from $\M_1^c(X)$ such that $\mu(A) =1$. Take an
automorphism $\g \in \hat U(R;\mu;1/2)\cap [T]$. Then there exists some $n$
such that $\mu(\{x\in X : Rx = T^nx\}) >0$. This contradicts the fact that
$\mu(E(R,T^n)) = 1\ \forall n\in \Z$. \hfill$\square$
\\

Given $T\in \Aut$, a Borel set $A\subset X$ is called a {\it complete
section} (or simply a {\it $T$-section}) if every $T$-orbit meets $A$ at
least once. If there exists a complete Borel section $A$ such that $A$
meets every $T$-orbit exactly once, then $T$ is called {\it smooth}. In
this case, $X = \bigcup_{i\in \Z} T^iA$ and all the $T^iA$'s are disjoint.
The set of smooth automorphisms is denoted by ${\cal S}m$.

A measurable set $W$ is said to be wandering with respect to $T\in \Aut$ if
the sets $T^nW,\ n\in\Z$, are pairwise disjoint. The $\sigma$-ideal
generated by all $T$-wandering sets in $\B$ is denoted by ${\mathcal
W}(T)$. By Poincar\'{e} recurrence lemma, one can state that given $T\in
\Aut$ and $A\in \B$ there exists $N\in {\mathcal W}(T)$ such that for each
$x\in A\setminus N$ the points $T^nx$ return to $A$ for infinitely many
positive $n$ and also for infinitely many negative $n$ [N]. The points from
the set $A\setminus N$ are called {\it recurrent}.

Assume that all points from a given set $A$ are recurrent for a Borel
automorphism $T$. Then for $x\in A$, let $n(x) = n_A(x)$ be the smallest
positive integer such that $T^{n(x)}x \in A$ and $T^ix \notin A,\ 0< i <
n(x)$. Let $C_k = \{x\in A\  \vert\  n_A(x) =k\},\ k\in \N $, then $ T^kC_k
\subset A $ and $ \{T^iC_k\  \vert\ i=0,...,k-1\} $ are pairwise disjoint.
Note that some $C_k$'s may be empty. Since $T^nx\in A$ for infinitely many
positive and negative $n$, we obtain
$$
 \bigcup_{n\ge0}T^nA =  \bigcup_{n\in \Z}T^nA = X
 $$
and
$$
 \bigcup_{n\ge0}T^nA  = \bigcup_{k=1}^\infty\bigcup _{i=0}^{k-1} T^iC_k.
$$
This union decomposes $X$ into $T$-towers $\xi_k = \{T^iC_k \ \vert \
i=0,..., k-1\},\ k \in \N$, where $C_k$ is the base and $T^{k-1}C_k$ is the
top of $\xi_k$. In particular, the number of these towers may be finite.

The next lemma is one of the main tools in our study of Borel automorphisms.

\begin{lemma}\label{L2.1} Let $T\in \Aut$ be an aperiodic Borel automorphism
of a standard Borel space $\bs$. Then there exists a sequence $(A_n)$ of
Borel sets such that\\ (i) $X=A_0 \supset A_1 \supset A_2\supset \cdots,$\\
(ii) $\bigcap_n A_n =\emptyset,$\\ (iii) $A_n$ and $X\setminus A_n$ are
complete $T$-sections, $n\in \N$,\\ (iv) for $n\in \N$, every point in
$A_n$ is recurrent,\\ (v) for $n\in \N$,\ $A_n\cap T^i(A_n) =\emptyset,\
i=1,...,n-1$,\\ (vi) for $n\in \N$, the base $C_k(n)$ of every non-empty
$T$-tower is an uncountable Borel set, $k\in \N$.
\end{lemma}

\noindent {\it Proof}. See [BeKe, Lemma 4.5.3] where (i) - (iii) have been
proved in more general settings of countable Borel equivalence relations. It
is shown in [N, Chapter 7] that one can refine the choice of $(A_n)$ to get
(iv) and (v). It is clear that one can remove an at most countable set of
points from each $A_n$ to prove (vi). \hfill{$\square$}
\\

\begin{definition}\label{markers} A sequence of Borel sets satisfying
conditions (i) - (vi) of Lemma $\ref{L2.1}$ is called a vanishing sequence
of markers.
\end{definition}

Note that usually $(A_n)$ is called a {\it vanishing sequence of markers}
if it satisfies (i) - (iv). We have added two more conditions, (v) and
(vi), which we will need in the constructions in Section 4.

\begin{remark}\label{change} {\rm We will use below the following
{\it changing-of-topology} result (see, for example, [Ke1, N]). Let $T\in
\Aut$ and let $(\xi_n)$ be a sequence of at most countable partitions of
$X$ such that: (a) $\xi_{n+1}$ refines $\xi_n$; (b) $\bigcup_n \xi_n$
generates the $\sigma$-algebra of Borel sets $\B$. Then we may introduce a
topology $\omega$ on $X$ such that: (i) $(X, \omega)$ becomes a Polish
0-dimensional space, (ii) $\B(\omega)= \B$ where $\B(\omega)$ is the
$\sigma$-algebra generated by $\omega$-open sets, (iii) all elements of
partitions $\xi_n, n \in \N$  are clopen in $\omega$, (iv) $T$ is  a
homeomorphism of $(X,\omega)$. In particular, by changing-of-topology, we
can choose the elements of the partitions corresponding to a vanishing
sequence of markers to be clopen.}
\end{remark}

\begin{proposition}\label{P2.2}  Let $T\in \Aut$ be an aperiodic Borel
automorphism of a standard Borel space $ \bs $. Then there exists a sequence
of periodic automorphisms $(P_n)$ of $ \bs $ such that $P_n
\stackrel{\tau}{\longrightarrow} T,\ n\to \infty$. Moreover, the $P_n$ can
all be taken from $[T]$.
\end{proposition}

\noindent {\it Proof}. If $T$ is a smooth automorphism, then the proof is
obvious. Let $(A_n)$ be a vanishing sequence of markers for $T$. Then, as we
have seen above, $A_n$ generates a decomposition of $X$ into $T$-towers
$\xi_k(n) = \{T^iC_k(n) \ \vert\ i=0,...,k-1\}$ and $\bigcup_k C_k(n) =
A_n$. Define
$$
P_nx =\left\{ \begin{array}{ll}
Tx, & {\rm if}\  x\notin B_n= \bigcup_{k=1}^\infty T^{k-1}C_k(n)\\
\\
T^{-k+1}x, & {\rm if}\  x\in T^{k-1}C_k(n),\ {\rm for\ some}\ k
\end{array} \right.
$$
Then $P_n$ belongs to $[T]$ and the period of $P_n$ on $\xi_k(n)$ is
$k$.  Note that $P_n$ equals $T$ everywhere on $X$ except $B_n$, the
union of the tops of the towers.

It follows from Lemma \ref{L2.1} that $(B_n)$ is a decreasing sequence of
Borel subsets such that $\bigcap_n B_n =\emptyset$. This means that for any
$x\in X$
there exists $n(x)$ such that $x\notin B_{n},\ n \ge n(x)$. Moreover, if for
some $x\in X$, $P_nx = Tx$, then $P_{n+1}x = Tx$. These facts prove that for
each $x$, all the $P_nx$ are eventually the same and equal to $Tx$, that is
$P_n$ converges to $T$ in $\tau$. \hfill{$\square$}
\\

We now give a version of the Rokhlin lemma for aperiodic Borel
automorphisms. We should also remark that B.~Weiss proved a measure-free
version of the Rokhlin lemma [W].

\begin{theorem}\label{T2.3} {\rm (Rokhlin lemma)}. Let $m\in \N$ and let $T$
be an aperiodic Borel automorphism of $\bs$. Then for any $\e>0$ and any
measures $\mu_1,...,\mu_p$ from $ {\cal M}_1(X)$ there exists a Borel subset
$F$ in $X$ such that $F, TF,..., T^{m-1}F$ are pairwise disjoint and
$$
\mu_i(F\cup TF\cup\cdots \cup T^{m-1}F) > 1-\e,\ \ i=1,...,p.
$$
\end{theorem}

\noindent {\it Proof}. We will use notation from the proof of Proposition
\ref{P2.2}. Clearly, it suffices to consider the case of non-smooth
automorphism $T$ only. Let $(A_n)$ be a vanishing sequence of markers. Note
that for any $\mu\in {\cal M}_1(X)$, $\mu(A_n) \to 0$ as $n\to \infty$
because $A_n$ decreases to the empty set. By the same reasoning,
$\mu(B_n)\to 0$. For every $n$, the space $X$ can be represented as a union
of $T$-towers $\xi_k(n)$ where the height of $\xi_k(n)$ is $k$ (see the
proof of Proposition \ref{P2.2}. Let
$$
D_n(m)= \bigcup_{k=1}^{m-1}\xi_k(n).
$$
Since $D_n(m)\subset \bigcup_{k=0}^{m-2}T^kA_n $, we see that there
exists $n_0$ such that  for $n>n_0$
\begin{equation}\label{2.1}
\mu_i(D_n(m)) < \frac{\e}{2},\ \ i=1,...,p.
\end{equation}
Let $B'_n= \bigcup_{k\ge m}T^{k-1}C_k(n)$. Similarly, we can deduce that for
all sufficiently large $n$
\begin{equation}\label{2.2}
\mu_i(B'_n\cup T^{-1}B'_n\cup \cdots \cup T^{-m+2}B'_n) \le \frac{\e}{2},\ \
i=1,...,p.
\end{equation}
Let $n$ be chosen so large that (\ref{2.1}) and (\ref{2.2}) hold. Define
$F$ by the
following rule. In each $T$-tower $\xi_k(n),\ k\ge m$, we take every $m$-th
set beginning with $C_k(n)$, i.e.
$$
 F = \bigcup_{k\ge m}\bigcup_{j=0}^{[\frac{k}{m}] -1}T^jC_k(n).
 $$
Then $F\cap T^jF =\emptyset,\ j=1,...,m-1$, and
$$
\mu_i(X - (F\cup TF\cup\cdots \cup T^{m-1}F)) < \e,\ \ i=1,...,p,
$$
in view of (\ref{2.1}) and (\ref{2.2}).\hfill{$\square$}
\\

It follows from Theorem \ref{T2.3} that in any $\tau$-neighborhood
$U(T;\mu_1,...,\mu_n;\e)$ of an aperiodic Borel automorphism $T$ there
exists a pointwise periodic automorphism. Thus we obtain the following
corollary from the above results:

\begin{corollary}\label{C2.4} $(1)$ The sets $\per$ and $\per_0$ are dense
in $(\Aut,\tau)$. Moreover, $\per$ is also dense in $\Aut$ with respect to
topologies $\tau'',\overline{p}$.\\ $(2)$ $\per\cap [T]$ is $\tau$-dense in
$[T]$ for each aperiodic $T$.\\ $(3)$ The set $\per$ is not dense in $\Aut$
with respect to $p$.
\end{corollary}
{\it Proof}. (1) and (2) are obvious. To prove (3), take an uncountable
Borel set $E$ and an aperiodic Borel automorphism $T$ such that
$TE\subsetneq E$. Then $W(T;E)$ has no periodic automorphisms.
\hfill$\square$
\\

We observe that the following result can also be proved. The details are
left to the reader.

\begin{corollary}\label{C2.5} Let $N\in \N$ and let a
$\tau$-neighborhood $U = U(T;\mu_1,...,\mu_n:\e)$ be given. Define $V = U(T;
(\nu^k_i)_{1\le i\le n,\ \vert k\vert \le N};\ \de)$ where $\nu^k_i =
\mu_i\circ T^{-k}$, $\de =(2N)^{-1}\e$. Then $V\subset U$ and for any $S\in
V$, we have that $S^j \in U(T^j; \mu_1,...,\mu_n,\e), \ j=1,...,N$.
\end{corollary}

We will need below the following statement proved in [BM].

\begin{lemma}\label{smooth} The set ${\cal S}m$ of smooth automorphisms is
dense in $\Aut$ with respect to the topology $p$.
\end{lemma}

Denote by $\ap \mod(Ctbl)$ the subset of $\Aut$ consisting of automorphisms
which are aperiodic everywhere except in an at most countable subset of
$X$.

\begin{theorem}\label{T2.6}  $(1)$
$\overline{\ap}^{\tau_0} = \ap \mod(Ctbl)$.\\ $(2)$ $\ap$ is a
nowhere dense closed subset in $(\Aut,\tau)$.
\end{theorem}

\noindent {\it Proof}. (1) We first show that the set $\ap \mod(Ctbl)$ is
closed with respect to $\tau_0$. Suppose $R\in \overline{(\ap
\mod(Ctbl))}^{\tau_0}\setminus (\ap \mod(Ctbl))$. Then there exist some
$m\in \N$ and an uncountable $R$-invariant Borel set $B$ such that $R$ has
period $m$ on $B$. Let $\mu$ be a continuous Borel probability measure such
that $\mu(B) =1$ and $\mu\circ R =\mu$. The $\tau_0$-neighborhood $U_0(R) =
U_0(R;\mu; \e/m)$ contains an automorphism $S$ from $\ap \mod(Ctbl)$. We
have that
$$
\mu(\{x\in B : Sx = Rx, ..., S^{m-1}x = R^{m-1}x, S^{m}x = R^{m}x =x\}) >
1-\e.
$$
In other words, $S$ is periodic on an uncountable Borel set, a contradiction.
Thus, we have that $\overline{\ap}^{\tau_0}\subset \ap \mod(Ctbl)$.

Conversely, if $R\in \ap \mod(Ctbl)$, then we need to show that any
$\tau_0$-neighborhood $U_0(R)$ of $R$ contains some aperiodic automorphism.
Indeed, the periodic part of $R$ is supported by an either countable or
finite set $A$. It is clear that if $A$ is infinite, then one can change $R$
on $A$ to produce an aperiodic automorphism from $U_0(R)$. If $A$ is finite,
then we take  a single aperiodic orbit $Orb_R(x),\ x\notin A$, and consider
the infinite set $A\cup Orb_R(x)$.

(2) It follows from Proposition \ref{P2.1}(3) and Corollary
\ref{C2.4} that $\Aut \setminus \ap$ is an open dense subset.
Therefore, $\ap$ is a closed nowhere dense set in $\tau$.
\hfill$\square$
\\

Let $(X, d)$ be a compact metric space (in particular, $X$ can be a Cantor
set). Recall that in this case we can define the metric $D$ on $\Aut$ as in
Remark \ref{rem}. We proved in [BDK] that the set of aperiodic
homeomorphisms is $D$-dense in $Homeo(X)$ when $X = \Om$ is a Cantor set.
Here we will find the closure of $\ap$ in the group $\Aut$ with respect to
the metric $D$.

Let $T$ be a Borel automorphism of $X$. Then $X$ is decomposed into the
canonical $T$-invariant partition $(Y_1,Y_2,...,Y_\infty)$ where $T$ has
period $n$ on $Y_n$ and $T$ is aperiodic on $Y_\infty$. We call $T$ {\it
regular} if all sets $Y_i,\ 1\leq i < \infty$, are uncountable.

\begin{lemma}\label{1case} Suppose that $T\in \Aut$ is regular. Then for any
$\e> 0$ there exists $S\in \ap$ such that $D(S,T) < \e$.
\end{lemma}
{\it Proof}. Since $TY_i = Y_i$, it suffices to find an aperiodic
automorphism $S$ satisfying the condition of the lemma for each set $Y_i$.
We can write down $Y_i$ as $E_i \cup TE_i\cup \cdots\cup T^{i-1}E_i$ for
some $E_i\in \B_0$. Let $(E_i(1),...,E_i(k_i))$ be a partition of $E_i$
into uncountable Borel sets such that ${\rm diam}(E_i(j)) < \e$ for all $j
=1,...,k_i$. Let $R(j)$ be an aperiodic automorphism of $E_i(j)$. Set for
$j=1,...,k_i$
$$
Sx = \left\{
\begin{array}{ll}
    Tx, & \mbox{$x\in \bigcup_{k=1}^{i-1} T^kE_i(j)$}\\
    \\
    TR(j)x, & \mbox{$x\in E_i(j).$}
\end{array}
\right.
$$
Clearly, $D(S,T) < \e$ on $Y_i$. \hfill$\square$

\begin{corollary} The set of aperiodic automorphisms from
$\widehat{Aut}(X,\B)$  is dense with respect to $\hat D$.
\end{corollary}

To answer the question when a non-regular automorphism $T\in \Aut$ belongs
to $\overline{\ap}^D$, we need to introduce the following definition. We
say that an automorphism $T$ is {\it semicontinuous} at $x\in X$ if for any
$\e >0$ there exists $z\neq x$ such that $d(x,z) < \e$ and $d(Tx,Tz) < \e$.

\begin{theorem}\label{D1} Let $T$ be a non-regular Borel automorphism from
$\Aut$ and let $Y_0$ denote the set $\bigcup_{i\in I}Y_i$ where each $Y_i$
is an at most countable set, $i\in I \subset \N$. Then $T\in
\overline{\ap}^D$ if and only if for every $x\in Y_0$ there exists $y\in
Orb_T(x)$ such that $T$ is semicontinuous at $y$.
\end{theorem}
{\it Proof}. We first suppose that for any $\e >0$ there exists an
aperiodic automorphism $S = S_\e$ such that $D(T,S) < \e$. Notice that the
fact that $T$ is pointwise periodic on $Y_0$ implies that $Orb_T(x) \neq
Orb_S(x)$ for any $x\in Y_0$. Hence there exists $y\in Orb_T(x)$ such that
$Sy\neq Ty$. On the other hand, we have that
\begin{equation}\label{ii}
d(Ty,Sy) <\e, \ \ \ \ d(T^{-1}(Sy), S^{-1}(Sy))< \e.
\end{equation}
Denoting by $z = T^{-1}Sy$, we obtain from (\ref{ii}) that $z\neq y$,
$d(z,y) < \e$, and $d(Tz,Ty) <\e$, that is $T$ is semicontinuous at
$y$.\footnote{We have not used in this part of the proof the fact that $T$
is non-regular.}

Suppose now that for every $x\in Y_0$ there exists $y\in Orb_T(x)$ such
that $T$ is semicontinuous at $y$. We need to show that for any $\e
>0$ there exists $S\in \ap$ such that $D(T,S) <\e$. It is clear that $S$
can be taken to coincide with $T$ on $Y_\infty$. Therefore, we need to
define $S$ on the at most countable set $Y_0$. We assume here that $Y_0$ is
infinite. It will be clear from the proof how one can deal with the case
when $Y_0$ is finite.

Take a finite partition $(C_1,...,C_n)$ of $X$ into Borel sets such that
${\rm diam}(C_i) < \e/2$ for all $i$. Denote by $A_{ij} = Y_0 \cap C_i \cap
T^{-1}C_j,\ i,j =1,...,n$. For each $x\in A_{ij}$, choose $y(x)\in
Orb_T(x)$ such that $T$ is semicontinuous at $y(x)$. The set $Y_0' =\{y(x)
: x\in Y_0\}$ is a subset of $Y_0$ intersecting each $T$-orbit of $x\in
Y_0$ exactly once. Set $A'_{ij} = A_{ij} \cap Y_0'$. Let $J =
\{(i_1,j_1),...,(i_p,j_p)\}$ be the set of those pairs $(i,j)$ for which
$A'_{ij} \neq \emptyset$. Then
\begin{equation}\label{Y}
Y_0 = \bigcup_{(i,j)\in J}\ \bigcup_{y(x)\in A'_{ij}}Orb_T(y(x)).
\end{equation}

Fix $(i,j) = (i_1,j_1)$. We have two possibilities: (a) $|A'_{ij}|
=\infty$, (b) $|A'_{ij}| < \infty$. If (a) holds, write down $A'_{ij}$ as
$\{...y_{-1},y_0,y_1,...\}$. Define $S$ on $\bigcup_{k\in \Z} Orb_T(y_k)$.
Set
$$
Sz = \left\{%
\begin{array}{ll}
    Tz, & \hbox{$z\in Orb_T(y_k),\ z\neq y_k$}\\
    Ty_{k+1}, & \hbox{$z = y_k$},\ \ \qquad k\in \Z.
\end{array}%
\right.
$$
In such a way, the set $\bigcup_{k\in \Z} Orb_T(y_k)$ is included in an
infinite $S$-orbit.

If (b) holds, then $A'_{ij} = \{z_1,...,z_q\}$. Let $\eta_1 = \min\{
d(z_i,z_j) : i\neq j,\ i,j = 1,...,q\}$ and let $0 <\eta < \min\{\e/2,
\eta_1\}$. By the hypothesis of the theorem, there exists $z\in C_i, z\neq
z_1,$ such that $Tz\in C_j$, $d(z,z_1) <\eta$, $d(Tz,Tz_1) < \eta$, and
the $T$-orbit of $z$ is infinite. To produce an $S$-orbit defined on
$\bigcup_{k=1}^q Orb_T(z_k)$, we can insert the $T$-orbits of
$z_1,...,z_q$ into $Orb_T(z)$. To do this, set $Sz = Tz_1$, $Sw = Tw, w\in
Orb_T(z), w\neq z$, and
$$
Sw = \left\{%
\begin{array}{ll}
    Tw, & \hbox{$w\in \bigcup_{k=1}^q Orb_T(z_k),\ w\neq z_1,...,z_q$}\\
    Tz_{i+1}, & \hbox{$w = z_i,\ i=1,...,q-1$} \\
    Tz , & \hbox{$w = z_q$.} \\
\end{array}%
\right.
$$
By the choice of the $C_i$'s, we see that in both cases (a), (b)
\begin{equation}\label{<}
d(Tx, Sx) + d(T^{-1}x, S^{-1}x) < \e
\end{equation}
on the set $\bigcup_{y\in A'_{ij}} Orb_T(y)$.

Take $(i_2,j_2)$ from $J$. By definition of $Y'_0$, we notice that
$\bigcup_{y\in A'_{ij}} Orb_T(y)$ does not meet the set of $T$-orbits going
through the points from $A'_{i_2,j_2}$. Therefore, we can apply
consequently the above procedure until the automorphism $S$ is defined
everywhere on $Y_0$. By (\ref{Y}) and (\ref{<}), we obtain that $D(T,S) <
\e$. \hfill$\square$
\\


{\bf 2.2. Incompressible automorphisms}. Let $T$ be an aperiodic Borel
automorphism of $\bs$. Let us denote by $[T]_0$ the set of Borel bijections
$\g: A\to B$ where $A, B$ are Borel subsets of $X$ and $\g x \in Orb_T(x),\
x\in A$. We call $A$ and $B$ equivalent with respect to $T$, $A\sim_T B$,
if there exists $\g\in [T]_0$ such that $\g(A) = B$. If there exists a
Borel subset $A$ of $X$ such that $X\sim_T A$ and $X\setminus A$ is a
complete $T$-section, then $T$ is called {\it compressible}. Otherwise, $T$
is called {\it incompressible}. We denote the set of incompressible
aperiodic automorphisms by $\inc$. It was proved in [DJK] that $T$ is
compressible if and only if $[T]$ contains a smooth aperiodic automorphism.
Let $M_1(T)$ denote the set of Borel probability $T$-invariant measures.
Clearly, for some automorphisms this set may be empty. For example, if $T$
is smooth, then $M_1(T) =\emptyset$. M.~Nadkarni [N] proved that $T$ is
incompressible if and only if there exists a $T$-invariant Borel
probability measure, i.e. $M_1(T) \neq \emptyset$.

\begin{theorem}\label{inc} The set $\inc$ is a closed nowhere dense subset
of $\ap$ with respect to the topology $p$.
\end{theorem}
{\it Proof}. We first show that $\inc$ is $p$-closed. Let $T\in
\overline{\inc}^p$. Choose a sequence $(\xi_n)$ of partitions of $X$ such
that: (i) $\xi_n = (F_n(1),...,F_n(k_n))$, $F_n(i) \in \B_0$; (ii)
$\xi_{n+1}$ refines $\xi_n$; (iii) $\bigcup_n \xi_n$ generates the
$\sigma$-algebra of Borel sets $\B$. By changing-of-topology results
(Remark \ref{change}), we may choose a topology $\omega$ on $X$ such that:
(i) $X$ is a Polish 0-dimensional space, (ii) $\B(\omega)= \B$ where
$\B(\omega)$ is the $\sigma$-algebra generated by $\omega$-open sets, (iii)
the sets $F_n(i)),\ n\in \N,\ i=0,1,...,k_n$ are clopen in $\omega$, (iv)
$T$ is a homeomorphism of $(X,\omega)$.

Denote by $W_n(T)$ the $p$-neighborhood $W(T; F_n(1),...,F_n(k_n))$. Then
$W_n(T)$ meets the set $\inc$ for every $n$. Let $S_n \in \inc \cap W_n(T)$
and let $\mu_n$ be an $S_n$-invariant probability measure. Then, for every
$n\in \N$,
\begin{equation}\label{1}
S_nF_n(i) = TF_n(i),\ \ i=1,...,k_n.
\end{equation}
Notice that if $m> n$, then  $F_n(i) = \bigcup_{j\in I} F_m(j)$ for each
$i= 1,...,k_n$ where $I \subset \{1,...,k_m\}$. Hence by (\ref{1})
$$
TF_n(i) = \bigcup_{j\in I} TF_m(j)=\bigcup_{j\in I} S_mF_m(i) = S_mF_n(i)
$$
and
\begin{equation}\label{2}
\mu_{m}(TF_n(i)) = \mu_m(F_n(i))\ \ \ m\geq n.
\end{equation}
The set $\{\mu_n : n\in \N\}$ contains a subsequence $(\mu_{n_k})$ which
converges to a Borel probability measure $\mu$ in the weak$^*$ topology.
Let us show that $\mu$ is $T$-invariant. For $B\in \bigcup_n \xi_n$ we
obtain that $\mu_{n_k}(B) \to \mu(B)$ and $\mu_{n_k}(TB) \to \mu(TB)$ as
$n_k\to \infty$, since $B$ and $TB$ are clopen sets (see e.g. [Bil]). It
follows from (\ref{2}) that for those sets $B$
$$
\mu(TB) = \lim_{n_k \to \infty} \mu_{n_k}(TB) = \lim_{n_k \to \infty}
\mu_{n_k}(B) =\mu(B).
$$
Since $\bigcup_n \xi_n$ generates $\B$, we see that $\mu$ is $T$-invariant.

To finish the proof, we refer to Lemma \ref{smooth} which provides us with
the following result: $\overline{{\cal S}m}^p = \Aut$. Since ${\cal S}m
\cap \inc = \emptyset$, we are done. \hfill$\square$
\\

Observe that in [BDK] we proved that if $X$ is a Cantor set, then
$\overline{Homeo(X)}^\tau = \Aut$. From Theorem \ref{inc}, we obtain that
$\overline{Homeo(X)}^p \subset \inc$.

Let $T\in \Aut$ be an aperiodic incompressible automorphism. Denote by
$$
{\rm Fix}(M_1(T)) = \{R\in \Aut : \mu\circ R = \mu, \ \forall \mu\in
M_1(T)\}
$$
and
$$
{\rm Pres}(M_1(T)) = \{R\in \Aut : R(M_1(T)) = M_1(T)\}.
$$

\begin{proposition}\label{fix} $(1)$ Let $T\in \inc$.
Then ${\rm Fix}(M_1(T))$ is closed in the topologies $\tau, \tau''$, and
$p$.\\
$(2)$ Let $T$ be an incompressible Borel automorphism of a compact metric
space $(X,d)$. Then ${\rm Fix}(M_1(T))$ is closed in the topology defined
by the metric $D$.\\
 $(3)$ The set ${\rm Pres}(M_1(T))$ is closed in $p$.
\end{proposition}
{\it Proof}. (1) We will first prove the statement for the topology
$\overline{\tau}$ which is equivalent to $\tau''$ by (\ref{tau''}). It will
follow from this result that  ${\rm Fix}(M_1(T))$ is closed in $\tau$. Let
$R\in \overline{{\rm Fix}(M_1(T))}^{\overline{\tau}}$ and $\mu \in M_1(T)$.
Then the $\overline{\tau}$-neighborhood $\overline{V}(R;\mu;\e) = \{S\in
\Aut : \sup_{F\in \B} |\mu(RF) - \mu(SF)| < \e\}$ meets ${\rm Fix}(M_1(T))$
for any $\e >0$. Therefore, for any Borel set $F$, we have that $|\mu (RF)
- \mu(F)| < \e$. Thus, $\mu\circ R = \mu$ and $R\in {\rm Fix}(M_1(T))$.

The fact that ${\rm Fix}(M_1(T))$ is closed in $p$ can be proved similarly
(see also (3)).

(2) We need to show that for every automorphism $R\in \overline{{\rm
Fix}(M_1(T))}^D$ and every $\mu \in M_1(T)$ one has $\mu\circ R = \mu$.
Take a sequence $(\g_n)$ from ${\rm Fix}(M_1(T))$ such that $D(\g_n, R) \to
0$ as $n\to\infty$. For $A\subset X$ and $\alpha > 0$, denote by
$B_\alpha(A)$ the $\alpha$-neighborhood of $A$. We notice that
$\overline{A}^d = \bigcap_{\alpha >0}B_\alpha(A)$ (a countable
intersection). It is clear, that for any $\alpha >0$, there exists
$n_\alpha$ such that $\g_n(A) \subset B_\alpha(RA)$ for all $n\geq
n_\alpha$.

Fix $\e >0$ and $\mu \in M_1(T)$. The following statement follows easily
from Luzin's theorem.
\medskip

\noindent {\bf Claim} Let $R,\mu, (\g_n)$ and $\e$ be as above. Then there
exists a closed subset $F_\e$ of $X$ such that the automorphisms $R$,
$R^{-1}$, and $\g_n$, being restricted to $F_\e$, are homeomorphisms and
$\mu(\tilde F_\e) > 1-\e$ where $\tilde F_\e = RF_\e\cap F_\e \cap
R^{-1}F_\e$.
\medskip

Clearly, one can choose the sets $\tilde F_\e$ such that $\tilde F_{\e_1}
\subset \tilde F_\e$ when $\e > \e_1$.  Let $C$ be a closed subset of
$\tilde F_\e$. Then for $\alpha
> 0$ and sufficiently large $n$, $\mu(C) = \mu(\g_nC) \leq
\mu(B_\alpha(RC))$. Hence
\begin{equation}\label{?}
\mu(C)\leq \lim_{\alpha \to 0}\mu(B_\alpha(RC)) = \mu (\bigcap_{\alpha >0}
B_\alpha(RC))= \mu(\overline{RC}^d) =\mu(RC).
\end{equation}
Similarly to (\ref{?}) we can show that $\mu(C) \leq \mu(R^{-1}C)$ and
therefore $\mu(A) = \mu(RA)$ for any Borel set $A \subset \tilde F_\e$.

Let now $\tilde F = \bigcup_{\e>0} \tilde F_\e$ (a countable union).
Clearly, $\mu(\tilde F) =1$. The above argument shows that $\mu(A) =
\mu(RA)$ for any Borel set $A \subset \tilde F$. It remains to check that
if $E$ is a Borel subset of $X\setminus \tilde F$, then $\mu(RE) = 0$.
Indeed, given $\alpha >0$, we can find a sufficiently large $n= n(\alpha)$
such that $RE \subset B_\alpha(\g_nE)$. Since $\lim_{\alpha\to 0}
\mu(B_\alpha(\g_nE))= 0$, we have that $\mu(RE) = 0$.

(3) Let us show that ${\rm Pres}(M_1(T))$ is closed in the topology $p$.
Indeed, if $R \in \overline{{\rm Pres}(M_1(T))}^p \setminus {\rm
Pres}(M_1(T))$, then there exist $\mu_0 \in M_1(T)$ and a Borel set $E$
such that $\mu_0(RTE) \neq \mu_0(RE)$. On the other hand, the
$p$-neighborhood $W(R;E,TE)$ contains some $S\in  {\rm Pres}(M_1(T))$ and
therefore $RE = SE,\ RTE = STE$. Then  $\mu_0(RTE) = (\mu_0\circ S)(TE) =
\mu_0\circ S(E) = \mu_0(RE)$, a contradiction.  \hfill$\square$

\begin{remark}\label{odom}
{\rm (1) For $T\in \inc$, the full group $[T]$ is a subset of ${\rm
Fix}(M_1(T))$ and the normalizer $N[T] = \{S\in \Aut : S[T]S^{-1} = [T]\}$
is a proper subset of ${\rm Pres}(M_1(T))$.  On the other hand, we know
that in Cantor dynamics the set ${\rm Fix}(M_1(T))$ is the closure of
$[T]$ in $D$ [GPS2]. In Borel dynamics the situation is different. We
first recall the definition of odometers.

Let $\{\la_t\}_{t=0}^\infty$ be a sequence of integers such that $\la_t\geq
2$. Denote by $p_{-1}=1,\ p_t=\la_0\la_1\cdots \la_t,\ t=0,1,...\ .$ Let
$X$ be the group of all $p_t$-adic numbers; then any element of $X$ can be
written as an infinite formal series:
$$
X=\{x=\sum_{i=0}^\infty x_ip_{i-1}\ \vert\ x_i\in (0,1,...,\la_i-1)\}.
$$
It is well known that $X$ is a compact metric abelian group endowed with
the metric $d(x,y) = (n+1)^{-1}$ where $n = \min\{ i : x_i \neq y_i\},\ x=
(x_i), y=(y_i)$. By definition, an odometer $T$ is the transformation
acting on $X$ as follows: $T x=x+1,\ x\in X$, where $1
=1p_{-1}+0p_0+0p_1+\cdots \in X$\footnote{More general, we call $S$ an
odometer if $S$ is Borel isomorphic to $T$}. From topological point of
view, $(X,T)$ is a strictly ergodic Cantor system and the set $M_1(T)$ is
a singleton. The orbit $Orb_T(0)$ is dense in $X$, that is every $b\in X$
can be approximated in $d$ by integer adic numbers. If $b\in X$, then $T_b
: x \mapsto x+b$ commutes with $T$. It is known that the topological
centralizer $C(T)$ coincides with $\{T_b : b\in X\}$. Since $d(x + b, x
+c)= d(b,c),\ b,c \in X$, it follows from Proposition \ref{fix} that
$$
C(T) \subset \overline{\{T^n : n\in \Z\}}^D \subset \overline{[T]}^D
\subset {\rm Fix}(M_1(T))
$$
and}
$$
C(T) \ \setminus\ \{T^n : n\in \Z\} \subset {\rm Fix}(M_1(T))\ \setminus\
\{T^n : n\in \Z\}.
$$
\end{remark}

\noindent {\bf 2.3. Borel automorphisms of rank 1}. Recall the definition
of {\it rank} 1 Borel automorphisms following [N].

Let $T\in \Aut$ be an aperiodic non-smooth automorphism and let $\xi =
(B_0,...,B_n)$ be a $T$-tower, that is all $B_i$'s are disjoint where $B_i =
T^iB_0,\ i=1,...,n$. Then $B_0$ and $B_n$ are called the base and and the
top of $\xi$.

Suppose a disjoint collection $\eta= \bigcup_{j\in J} \xi(j)$ (finite or
countable) of $T$-towers is given where $\xi(j) = (B_0(j),...,B_{n_j}(j))$
and $B_0(j)\in {\cal B}_0$ is a Borel uncountable set,  $j\in J$. Then
$\eta$ is called a $T$-multitower and $Y= \bigcup_{j\in
J}\bigcup_{i=0}^{n_j} B_i(j)$ is called the support of $\eta$. The
cardinality of $J$ is called the the rank of the multitower. A multitower
$\eta'$ is said to refine $\eta$ if every atom $B'_{i'}(j')$ of $\eta'$ is
a subset of some atom $B_i(j)$ of $\eta$.

\begin{definition}\label{D2.7}  We say that $T$ has
{\rm rank at most $r$} if there exists a sequence $(\eta_n)$ of
$T$-multitowers of rank $r$ or less such that $\eta_{n+1}$ refines $\eta_n$
and the collection of all atoms in $\eta_n$, taken over all $n\in \N$,
generates $\B$. Then $Y_n\subset Y_{n+1}$ and $\bigcup_n Y_n = X$ where
$Y_n$ is the support of $\eta_n, \ n\in \N$. We say that $T$ has {\rm rank
$r$} if $T$ has rank at most $r$ but does not have rank at most $r-1$. If
$T$ does not have rank $r$ for any finite $r$ then, by definition, $T$ has
infinite rank. Denote by ${\cal R}(n)$ the set of automorphisms of rank $n$.
\end{definition}

A complete description of the structure of Borel automorphisms of rank 1
can be found in [N]. Here, we observe only that any $T\in {\cal R}(1)$ can
be obtained as a $\tau$-limit of a sequence $(T_n)$ of partially defined
Borel automorphisms. For this, we use the cutting and stacking method to
produce a refining sequence $(\xi_n)$ of towers satisfying Definition
\ref{D2.7}. More precisely, $\xi_n$ is first cut into $T_n$-subtowers
$\xi_n(k) = (C_n^0(k), \dots, C^{h_n}_n(k)),\ k=0,\dots, p_n,$ with
$T_n^i(C_n^0(k))= C_n^i(k),\ i= 1,\dots, h_n$, and then some spacers
$D_n^{h_n+1}(k),\dots, D_n^{m_n(k)}(k)$ are added to each subtower to
extend $\xi_n(k)$ to $\xi'_n(k)$. One defines $T_{n+1}(C^{h_n}_n(k)) =
D_n^{h_n+1}(k)$ and $T^i_{n+1}(D_n^{h_n+1}(k)) = D_n^{h_n+i+1}(k)),\ i=
1,\dots, m_n(k) - h_n- 1$. To construct $\xi_{n+1}$ one takes successively
the extended $T_{n+1}$-subtowers $\xi'_n(0),\dots, \xi'_n(p_n)$ and then
makes from them a single $T_{n+1}$-tower by concatenating those subtowers
and setting $T_{n+1}(D_n^{m_n(k)}(k)) = C_n^0(k+1),\ k=0,\dots,p_n$. Thus,
the base and the top of $\xi_{n+1}$ are $C_n^0(0)$ and
$D_n^{m_n(p_n)}(p_n)$ respectively. Remark that the spacers that enlarge
each $\xi_n(k)$ are taken from $X\setminus Y_n$. Finally, as $n\to\infty$,
we get a Borel automorphism $T$ of rank 1 as the limit of $T_n$.

If we assume that for given $T\in {\cal R}(1)$, from the above
construction, one can choose a sequence $(\xi_n)$ such that $Y_n = X$ for
every $n$ (no spacers can be added), then we get that $T$ belongs to the
simplest subclass of rank 1 Borel automorphism, the so called {\it
odometers}. The exact description of odometers is given in Remark
\ref{odom}. Let us denote this subclass by ${\cal O}d$.

We will need the following simple fact. If $T\in \Aut$ and $(F_1,...,F_n)$
is a Borel partition of $X$ such that $TF_1 = F_2,..., TF_{n-1} = F_n, TF_n
= F_1$, then the $p$-neighborhood $W(T; F_1,...,F_n)$ contains an odometer
$S$.

The next proposition shows that any rank 1 automorphism is a limit of
odometers in $\tau$. Later, in Section 4, this result will be strengthen.

\begin{theorem}\label{T2.8}   $\overline{{\cal R}(1)}^\tau = \overline{{\cal
O}d}^\tau$.
\end{theorem}

\noindent
{\it Proof}.  We need to show only that given $\e>0$, $T\in  {\cal R}(1)$,
and $\mu_1,...,\mu_n \in {\cal M}_1(X)$, there exists an odometer $S$ such
that $\mu_i(E(S,T))<\e$ for all $i$.

We first assume that every measure $\mu_i$ is continuous. Let $(\xi_n),\
\xi_n = (C_n^0, \dots, C^{h_n}_n)$, be a refining sequence of $T$-towers as
above. Then we see that $\mu_i(Y_n)\to 0$ as $n\to \infty$. Since atoms from
$(\xi_n)$ generate $\B$, we have that $\mu_i(C_n^0 \cup C^{h_n}_n) \to 0$ as
$n\to \infty$ for $i=1,...,n$. Find $N\in \N$ such that $\mu_i(Y_n \cup
C_n^0 \cup C^{h_n}_n) < \e$ for $n\ge N$. Clearly, we can define a new Borel
automorphism $S$ such that $Sx = Tx$ if $x\in \bigcup_{i=0}^{h_n-1} C^i_n$
and $S(C^{h_n}_n) = Y_n$, $S(Y_n) = C_n^0$. In other words, we have
constructed the $S$-tower $(C_n^0, \dots, C^{h_n}_n, Y_n)$ which partitions
$X$. Clearly, the definition of $S$ can be extended to produce an odometer
on $X$ which belongs to $U(T;\mu_1,...,\mu_n; \e)$.

Now suppose that every measure $\mu_i$ can have points of positive measure,
say $\{x_k(i)\}_{k\in \N},\ i=1,...,n$. Then find a finite set $Y = \{x_k(i)
: k\in I(\mu_i)\subset \N\}$ where a finite subset $I(\mu_i)$ is determined
by the condition
$$
\sum_{k\notin I(\mu_i)} \mu_i(\{x_k(i)\}) \le \frac{\e}{2},\ \
i=1,..., m.
$$
Notice that there exists a refining sequence of $T$-towers  $(\xi_n)$ such
that, for sufficiently large $n$, points from $Y$ do not belong to  the
base and the top of $\xi_n$. Indeed, it follows from the fact that
$(\xi_n)$ generates $\B$. The rest of the proof is the same as for
continuous measures.\hfill$\square$


\sect{Comparison of the topologies}.

\setcounter{equation}{0}

In this section, our main aim is to prove Theorem \ref{T1.2}, which
clarifies relationships between the topologies $\tau, \tau',\tau'',
\tau_0,\ p,\ p_0,\ \tilde p$, and $\overline{p}$.

\begin{proposition}\label{P3.1} $(1)$ The topology $\tau$ is strictly stronger
than $\tau_0$.\\ $(2)$ The topologies $p$ and $p_0$ are equivalent.
\end{proposition}
{\it Proof}. (1) We need to show only that $\tau_0$ is not equivalent to
$\tau$. For $y\in X$, take $U= U({\mathbb I}; \delta_y; 1/2)$. We will show
that for any $U_0 = U_0({\mathbb I}; \nu_1,...,\nu_n;\e)\ (\nu_i \in
\M_1^c(X))$ there exists a Borel automorphism $S$ from $U_0$ such that
$S\notin U$. To see this, take an uncountable Borel set $B$ such that $y\in
B$ and $\nu_i(B) < \e$ for all $i$. Let $S$ be a freely acting automorphism
on $B$ such that $Sx = x$ for $x\in B^c$.

(2) We will prove that any $p$-neighborhood $W(T;F_1,...,F_n)$ contains a
$p_0$-neighborhood $W_0(T;C_1,...,C_n)$ where the $C_i$'s are uncountable.
Without loss of generality, we can assume that $(F_1,...,F_n)$ is a
partition of $X$. Suppose $F_1,...,F_k$ are uncountable sets and
$F_{k+1},...,F_n$ are countable (or finite) ones. Define $C_i = F_i,\
i=1,...,k$, and $C_i = F_1 \cup F_i,\ i=k+1,...,n$. Clearly, if $R$ is a
Borel automorphism such that $RC_i = TC_i$ then $RF_i = TF_i$ for $i=
1,..., n$. \hfill$\square$

\begin{theorem}\label{T3.1}  The topologies $\tau $ and $\tau'$ are
equivalent.
\end{theorem}

\noindent {\it Proof}. It is sufficient to consider neighborhoods of
${\mathbb I}$ only since $\Aut$ is a topological group. Take a neighborhood
$U=U({\mathbb I};\ \mu_1,...,\mu_n;\ \e )$. We will show that
$U':=U'({\mathbb I};\ \mu_1,...,\mu_n;\ \e /4)\subset U$. By definition, if
$T\in U'$, then for any Borel set $F$ we have
\begin{equation}\label{3.1}
\mu_i(F\ \De \ T(F))<\e /4,\ \ \ \ i=1,...,n.
\end{equation}
We need to estimate $\mu_i(E(T,{\mathbb I}))$. Note that $ E(T,{\mathbb I})$
can be partitioned as a disjoint union $ X_2\cup X_3\cup \cdots \cup
X_\infty $ where the period of $T$ on $ X_k $ is $ k $ and $T$ is aperiodic
on $X_\infty$. Apply Theorem \ref{T2.3}
 with $X = X_\infty,\ m=2,$ and $\e/4$. We
obtain a Borel subset $F' \subset X_\infty$ such that $F'\cap TF'=\emptyset$
and by (\ref{3.1})
\begin{equation}\label{3.2}
 \mu_i(X_\infty)  <  \mu_i(F'\cup TF') + \frac{\e}{4}
  =  \mu_i(F'\ \De\ TF') +  \frac{\e}{4}  <  \dfrac{\e}{2}.
\end{equation}
For every $2\le k <\infty$, let us take $Y_k$ such that
$$
 X_k = \bigcup_{j=0}^{k-1} T^jY_k.
$$
Let $F_1 = \bigcup_{2\le k<\infty}Y_k$. Then $F_1\ \De\ TF_1= F_1\cup TF_1$
and therefore
\begin{equation}\label{3.3}
\mu_i(F_1) \le \mu_i (F_1\ \De\ TF_1) < \e/4,\ i=1,...,n.
\end{equation}
Denote by
$$
 F_2 = \bigcup_{2\le k<\infty}\ \bigcup_{0 \le j\le [\frac{k-1}{2}]} T^{2j}Y_k.
$$
Then $F_2\cap TF_2 \subset F_1$ and $ F_2\cup TF_2 = \bigcup_{2\le k<
\infty}X_k. $ Thus, we get from (\ref{3.1}) and (\ref{3.3})
 $$
 \mu_i(\bigcup_{2\le k< \infty}X_k) = \mu_i(F_2\ \De\ TF_2) + \mu_i(F_2\cap
TF_2) < \frac{\e}{2},\ i=1,...,n
$$
This result together with (\ref{3.2}) shows that $T\in U$ and therefore
$U'\subset U$.

Conversely, suppose $U'=U'({\mathbb I};\ \mu_1,...,\mu_n;\ \e )$ is
given. We will show that $U({\mathbb I};\ \mu_1,...,\mu_n;\
\e/2)\subset U'$. Indeed,  let $S\in  U({\mathbb I};\
\mu_1,...,\mu_n;\e/2 )$, then $\mu_i(E(S,{\mathbb I}))<\e/2 ,\
i=1,...,n$. Thus, for a Borel subset $F\subset E(S,{\mathbb I})$, we
have
 \begin{eqnarray*}
\mu_i(F\ \De\ S(F))& \le & \mu_i(F\cup S(F))\\
&\le & \mu_i(F)+\mu_i(SF)\\
& \le & 2\mu_i(E(S,{\mathbb I}))\\
& \le &\e,\ \ \ \ i=1,...n.
\end{eqnarray*}
If $F\subset X -E(S,{\mathbb I})$, then $F\ \De\ SF=\emptyset$.
Thus $U({\mathbb I};\ \mu_1,...,\mu_n;\ \e /2)\subset U'$. \hfill${\square}$

\begin{proposition}\label{P3.2} The topology $\tau$ (and therefore
$\tau'$) is strictly stronger than $\overline{p}$.
\end{proposition}

\noindent {\it Proof}. We define another topology, denoted $\tau'_0$, on
$\Aut$. By definition, $\tau'_0$ is generated by the base of neighborhoods
$$
U'_0(T; \mu_1,...,\mu_n; \e )= \{ S\in \Aut \ |\ \sup_{F\in {\cal
B}}\mu_i(TF\ \De \ SF)
$$
$$
 + \sup_{F\in {\cal B}}\mu_i(T^{-1}F\ \De \ S^{-1}F)<
\e,\ i=1,...,n\}.
$$
Obviously, $\tau'_0$ is stronger than $\tau'$.
\smallskip

\noindent {\bf Claim 1} {\it $\tau'_0$ is equivalent to $\tau$}.
\smallskip

In fact, we need to show only that $\tau$ is stronger than $\tau'_0$. This
assertion can be proved in the same method which was used to establish
that $\tau$ is stronger than $\tau'$ in Theorem \ref{T3.1}. Using this
fact, we obtain
$$
\tau \succ\tau'_0\succ\tau' \sim \tau
$$
and the claim is proved.

To finish this part of the proof, we note that $\tau'_0$ is clearly
stronger than $\overline{p}$.

Now we show that $\tau$ is strictly stronger that $\overline{p}$. For this,
we need to find a $\tau$-neighborhood $U$ of the identity such that for any
$\overline{p}$-neighborhood $ \overline{W} =\overline{W}({\mathbb I};
(F_i); (\mu_j);\e)$ there exists $S\in \overline{W}$ which is not in $U$.
Take $U= U({\mathbb I}; \de_{x_0}; 1/2)$. Then $S \notin U$ if and only if
$Sx_0 \neq x_0$. Thus, we have to show that in every $ \overline{W}$ there
exists $S$ such that $Sx_0\neq x_0$.
\smallskip

\noindent{\bf Claim 2} {\it Every $\overline{p}$-neighborhood of the
identity $\overline{W}({\mathbb I}; (F_i); (\mu_j);\e)$ contains a free
automorphism $S$}.
\smallskip

Indeed, if $S\in \overline{W}$ then $\mu_j(SF_i\ \De\ F_i) +
\mu_j(S^{-1}F_i\ \De\ F_i)< \e,\ i=1,...,n, j=1,...m$. Given $(F_i)$ and
$(\mu_j)$, one can find a freely acting $S$ satisfying the above condition.
To see this, we can assume that $X = {\mathbb R}$ and then $S$ can be taken
as a translation $x \to x+\alpha, \ \alpha \in {\mathbb R}$. The details are
left to the reader. \hfill$\square$
\\

Next, we will compare $\tau$ and $\tau''$. Our goal is to prove that the
uniform topology $\tau$ is strictly stronger than $\tau''$. To do this, we
will need a more convenient description of $\tau''$.

\begin{definition}\label{D3.3}  $(1)$ Let $\tilde\tau$ be the topology on
$\Aut$ defined by the base
$$
\widetilde V(T;\ \mu_1,...,\mu_n;\ \e )
$$
\begin{equation}\label{3.4}
=  \{S\in \Aut\ \vert\ \sup_{{\cal Q}}(\sum_{i\in
I}\vert\mu_j(TE_i) - \mu_j(SE_i)\vert) <\e,\ j=1,...,n\}
\end{equation}
where $\mu_1,...,\mu_n\in {\cal M}_1(X)$ and  supremum is taken over
all finite Borel partitions ${\cal Q} = (E_i)_{i\in  I}$ of $X$.

$(2)$ Define a new topology $\overline{\tau}$ on $\Aut$ using as
neighborhood base the sets
$$
\overline V(T;\ \mu_1,...,\mu_n;\ \e )
$$
\begin{equation}\label{3.5}
=  \{S\in \Aut\ \vert\ \sup_{F\in {\cal B}}\vert\mu_j(TF) -
\mu_j(SF)\vert <\e,\ j=1,...,n\}
\end{equation}
where $\mu_1,...,\mu_n\in {\cal M}_1(X)$.
\end{definition}

\begin{proposition}\label{P3.4} The three topologies $\tau'',
\tilde\tau$ and $\overline{\tau}$ are pairwise equivalent.
\end{proposition}

\noindent {\it Proof}. ($\tau'' \iff \tilde \tau$) Let $U'' = U''({\mathbb
I}; \mu_1,...,\mu_n; \e )$ be as in (\ref{1.3}). We will show that
$\widetilde V=\widetilde V ({\mathbb I}; \mu_1,...,\mu_n; \e /3) \subset
U''$ where $\widetilde V$ is defined in (\ref{3.4}). Take $S\in \widetilde
V$. For any $f\in B(X)_1,$ find $g(x)=\sum_{i\in I}a_i\chi_{E_i}$ such that
$$
\| f -g \| =\sup_{x\in X}\vert f(x) -g(x)\vert <\frac{\e}{3}.
$$
Here $\vert a_i\vert \le 1$ and $(E_i)_{i\in I}$ forms a partition of
$X$ into Borel sets.  Then for  $j=1,...,n$,
\begin{eqnarray*}
\left| \int_X (f(S^{-1}x)-f(x))d\mu_j\right| & \le & \left| \int_X
(f(S^{-1}x)-g(S^{-1}x))d\mu_j \right| \\
\\
& + & \left| \int_X (g(S^{-1}x)-g(x))d\mu_j\right|+\left| \int_X
(g(x)-f(x))d\mu_j\right|\\
\\
& \le & \frac{2\e }{3}+\left| \int_X (g(S^{-1}x)-g(x))d\mu_j\right|\\
\\
& = & \frac{2\e }{3}+\vert \sum_{i\in I}a_i(\mu_j(SE_i)-\mu_j(E_i))\vert\\
\\
& \le & \frac{2\e }{3}+\sum_{i\in I}\vert \mu_j(SE_i)-\mu_j(E_i)\vert
\\
& \le & \e .
\end{eqnarray*}
This proves that $\widetilde V \subset U''$.

Conversely, let $\widetilde V=\widetilde
V({\mathbb I};\mu_1,...,\mu_n;\e)$ be
given. It is sufficient to show that $U'' = U''({\mathbb
I};\mu_1,...,\mu_n;\e) \subset \widetilde V$. If $S\in U''$, then
for all $f\in B(X)_1$, we have
\begin{equation}\label{3.6}
\left| \int (f(S^{-1}x)-f(x))d\mu_j\right| <\e,\ \ j=1,...,n.
\end{equation}
We thus need to show that for every finite Borel partition ${\cal Q}=
(E_i)_{i\in I}$ and all $j=1,...,n$,
\begin{equation}\label{3.7}
\sum_{i\in I}\vert \mu_j(SE_i)-\mu_j(E_i)\vert <\e.
\end{equation}
For this, we consider the ${\cal Q}$-measurable function $g_j(x) =\sum_{i\in
I}a_i(j)\chi_{E_i}\in B(X)_1$, where for a finite set of measures $(\mu_j)$
and $S$ we define
$$
a_i(j)=\left\{ \begin{array}{rrr} 1, & {\mbox if}\ \mu_j(SE_i) >
\mu_j(E_i)\\ -1, & {\mbox if}\ \mu_j(SE_i) < \mu_j(E_i)\\ 0, & {\mbox if}\
\mu_j(SE_i) = \mu_j(E_i)
\end{array} \right.
$$
The definition of $g_j(x)$ and (\ref{3.6}) imply that for every $j=1,...,n$
\begin{eqnarray*}
\sum_{i\in I}\vert \mu_j(SE_i)-\mu_j(E_i)\vert & = & \left| \sum_{i\in
I}a_i(j)(\mu_j(SE_i)-\mu_j(E_i))\right| \\
\\
& = & \left| \int_X \sum_{i\in I}a_i(j)(\chi_{SE_i}-\chi_{E_i})d\mu_j\right|
\\
\\
& = & \left| \int_X (g_j(S^{-1}x)-g_j(x))d\mu_j\right|\\
\\
& < & \e.
\end{eqnarray*}
Thus, (\ref{3.7}) holds and therefore $U''\subset \widetilde V$.
\\

$(\tilde \tau \iff \overline{\tau})$ Let $\overline V ({\mathbb I};\
\mu_1,...,\mu_n;\ \e)$ be given as in (\ref{3.5}). Then, we will prove that
\begin{equation}\label{3.8}
\widetilde V ({\mathbb I};\ \mu_1,...,\mu_n;\ 2\e )\subset \overline
V ({\mathbb I};\ \mu_1,...,\mu_n;\ \e ).
\end{equation}
In fact, if $S\in \widetilde V ({\mathbb I};\ \mu_1,...,\mu_n;\ 2\e )$, then
for every  finite Borel partition ${\cal Q}$ into sets $(E_i)_{i\in I}$ we
have
$$
\sum_{i\in I}\vert \mu_j(SE_i)-\mu_j(E_i)\vert <2\e.
$$
In particular, if $F$ is a Borel set, then for ${\cal Q}=\{F, X -F\}$ and
$j=1,...,n$,
\begin{eqnarray*}
2\e & > & \vert \mu_j(SF)-\mu_j(F)\vert +\vert \mu_j(S(X -F)-\mu_j(X
-F)\vert\\
\\
& = & 2\vert \mu_j(SF)-\mu_j(F)\vert .
\end{eqnarray*}
Hence $\vert \mu_j(SF)-\mu_j(F)\vert <\e$ for any Borel set $F$ and
$j=1,...,n$, that is $S\in \overline V ({\mathbb I};\ \mu_1,...,\mu_n;\ \e
)$ and (\ref{3.8}) is proved.

Conversely, let $\widetilde V =\widetilde V({\mathbb I}; \mu_1,...,\mu_n;\
\e )$ be given as in (\ref{3.4}). We will prove that $\overline V ({\mathbb
I};\ \mu_1,...,\mu_n;\ \e /2)\subset \widetilde V$. Take $S\in \overline V
({\mathbb I};\ \mu_1,...,\mu_n;\ \e /2 )$. This means that for every $F\in
{\cal B}(X )$ and all $j=1,...,n$,
\begin{equation}\label{3.9}
\vert \mu_j(SF)-\mu_j(F)\vert <\e/2.
\end{equation}
Let ${\cal Q}=(E_i)_{i\in I}$ be a finite partition into Borel sets and
denote by $F_+(j)=\bigcup_{i\in I_+(j)}E_i$ where $I_+(j) =\{ i\in I\
\vert \ \mu_j (SE_i)\ge \mu_j (E_i)\}$. Then $F_-(j):= X -F_+(j)=
\bigcup_{i\in I_-(j)}E_i$, where $I_-(j)=\{ i\in I\ \vert \
\mu_j(SE_i)<\mu_j(E_i)\}$. We know, by (\ref{3.9}),
 that for $j=1,...,n$,
$$
\vert \mu_j(SF_+(j))-\mu_j(F_+(j))\vert <\e /2,
$$
$$
\vert \mu_j(SF_-(j))-\mu_j(F_-(j))\vert <\e /2.
$$
In other words,
\begin{eqnarray*}
\frac{\e}{2} & > & \vert \mu_j(SF_+(j))-\mu_j(F_+(j))\vert \\
\\
& = & \vert \sum_{i\in I_+(j)}(\mu_j(SE_i)-\mu_j(E_i))\vert \\
\\
& = & \sum_{i\in I_+(j)}\vert \mu_j(SE_i)-\mu_j(E_i)\vert .
\end{eqnarray*}
Similarly,
$$
\frac{\e}{2} >\vert \mu_j(SF_-(j))-\mu_j(F_-(j))\vert =\sum_{i\in
I_-(j)}\vert \mu_j(SE_i)-\mu_j(E_i)\vert .
$$
Therefore
$$
\sum_{i\in I}\vert \mu_j(SE_i)-\mu_j(E_i)\vert <\e,
$$
i.e. $S\in \widetilde V ({\mathbb I};\ \mu_1,...\mu_n;\ \e)$.
\hfill${\square}$

\begin{theorem}\label{T3.5} The topology $\tau$ (and therefore
$\tau'$) is strictly stronger than $\tau''$.
\end{theorem}

\noindent
{\it Proof}. The theorem will be proved in two steps. We first show that
$\tau$ is stronger than $\tau''$ (we give two proofs of this fact). Then
we prove that $\tau''$ cannot be equivalent to $\tau$.
\\

$1^0$. ($\tau \succ \tau''$) $1^{\rm st}$ proof. As mentioned above,
it suffices to consider neighborhoods of ${\mathbb I}$ only. Take
$U''=U''({\mathbb I};\ \mu_1,...,\mu_n;\ \e )$ as in (\ref{1.3}). We will
show that $U=U({\mathbb I};\ \mu_1,...,\mu_n;\ \e /2)\subset U''$.

Indeed, for $T\in U$, one has $\mu_i(E(T,{\mathbb I}))<\e /2,\ i=1,...,n$.
Then, for $f\in B(X)_1$,
\begin{eqnarray*}
\sup_{\Vert f\Vert \le 1}\vert \mu_i\circ T(f)-\mu_i(f)\vert & = & \sup_{\Vert
f\Vert \le 1}\left| \int_{X}(f(T^{-1}x)-f(x))d\mu_i\right| \\
\\
& = & \sup_{\Vert f\Vert \le 1}\left| \int_{E(T,{\mathbb
I})}(f(T^{-1}x)-f(x))d\mu_i\right| \\
\\
& \le & \sup_{\Vert f\Vert \le 1}2\Vert f\Vert \mu_i(E(T,{\mathbb
I}))<\e. \end{eqnarray*}
Thus, $U\subset U''$.
\\

($\tau \succ \tau''$) $2^{\rm nd}$ proof. By Theorem \ref{T3.1} and
Proposition \ref{P3.4}, the statement will be proved if we show that
$\tau'$ is stronger than $\overline{\tau} \sim \tau''$. To this end, we
note that for given $\mu_1,...,\mu_n\in {\cal M}_1(X)$ and $\e>0$, one has
\begin{equation}\label{3.10}
\overline V
({\mathbb I};\ \mu_1,...,\mu_n;\ \e )\supset U'({\mathbb I};\
\mu_1,...,\mu_n;\ \e ).
\end{equation}
Indeed, (\ref{3.10}) follows from (\ref{1.2}) and (\ref{3.5})
in view of the following
simple observation. If $S\in U'({\mathbb I};\ \mu_1,...,\mu_n;\ \e )$,
then
$$
\sup_{F\in {\cal B}}\mu_i(F\ \De\ SF)<\e,\ \ i=1,...,n.
$$
Since $\vert \mu_i(F)-\mu_i(SF)\vert \le\mu_i(F\ \De\ SF)$, we get that
$$
\sup_{F\in {\cal B}}\vert \mu_i(F)-\mu_i(SF)\vert <\e,
$$
i.e. $S\in \overline V({\mathbb I};\ \mu_1,...,\mu_n;\ \e)$.
\\

$2^0$. $(\tau \nsim \tau'')$ The theorem would be proved if we could show
that the following assertion holds. To see that $\tau'' \sim
\overline{\tau}$ cannot be equivalent to $\tau$, we should exhibit a
$\tau$-neighborhood, say $U({\mathbb I};\ \mu;\ \e)$, that does not contain
a $\overline{\tau}$-neighborhood. This means that we need to prove the
following claim.
\\

\noindent {\bf Claim 1} {\it There exists a $\tau$-neighborhood $U({\mathbb
I};\ \mu;\ \e)$ such that for every $\overline{\tau}$-neighborhood
$\overline V({\mathbb I};\ \nu_1,...,\nu_n;\ \de)=\overline V$ one can find
a Borel automorphism $S$ that belongs to $\overline V$ but does not belong
to $U({\mathbb I};\ \mu ;\ \e)$.}
\\

$2^0({\rm a})$. We first discuss the case where $\mu$ is a purely atomic
measure. The next claim shows that it is impossible to distinguish the
topologies $\tau$ and $\tau'' \sim \overline{\tau}$ by using
atomic measures only. \\

\noindent
{\bf Claim 2} {\it Let $U=U({\mathbb I};\ \mu_1,...,\mu_n;\ \e)$ be such
that all of the measures $\mu_i,\ i=1,...,n,$ are purely atomic. Then there
exists $\overline V= \overline V({\mathbb I};\ \nu_1,...,\nu_n;\ \de)$
such that $\overline V\subset U$.}
\\

To prove this claim,  fix a $\tau$-neighborhood $U({\mathbb I};\
\mu_1,...,\mu_n;\ \e)$  with atomic measures $\mu_i$. Then
there exists a set (at most countable) $\{ x^i_j\ |\ j\in \N ;\
i=1,...,n\}$ of points in $X$ such that $\mu_i(\{x^i_j\})>0$ and
$\sum_j\mu_i(\{x^i_j\})=1$. Choose $n_0$ such that for all $i=1,...,n,$
$$
\sum_{j>n_0}\mu_i(\{x^i_j\})<\e.
$$
For $\{ x^i_j\ \vert\ 1\le j\le n_0,\ 1\le i\le n\}$, define atomic measures
$\nu_i,\ i=1,...,n$:
$$
\nu_i(\{ x^i_j\} )=b^i_j>0
$$
where $\sum_{j=1}^{n_0}b^i_j=1$ and $b^i_j\ne b^i_{j_1}$ for all $i$
and $j\neq j_1$. Let
\begin{equation}\label{3.11}
\de <\min_{1\le i\le n}\ \ [\min_{j\ne j_1, 1\le j,j_1\le n_0}(\vert
b^i_j-b^i_{j_1}\vert,\ b^i_j)].
\end{equation}
Then we have
$$
\overline V({\mathbb I};\ \nu_1,...,\nu_n;\ \de)\subset U({\mathbb I};\
\mu_1,...,\mu_n;\ \e).
$$
Indeed, if $S\in \overline V({\mathbb I};\ \nu_1,...,\nu_n;\ \de)$, then
for each point $\{x^i_j\}\ (1\le j\le n_0, i=1,...,n)$ we get
\begin{equation}\label{3.12}
\vert \nu_i(\{x^i_j\})- \nu_i(\{Sx^i_j\})\vert <\de.
\end{equation}
Then (\ref{3.11}) and (\ref{3.12}) imply that $Sx^i_j= x^i_j$ for all $i$
and $1\le j\le
n_0$. In other words, $\nu_i$ is $S$-invariant. Then $x^i_j\notin
E(S,{\mathbb I})$ for all $i=1,...,n,\ j=1,...,n_0$. Therefore
$$
\mu_i(E(S,{\mathbb I})) \le \sum_{j>n_0} \mu_i(\{x^i_j\}) < \e,
$$
that is $S\in U({\mathbb I};\ \mu_1,...,\mu_n;\ \e)$. The claim is
proved.
\\

$2^0$(b). In view of Claim 2, we must consider continuous measures from
${\cal M}_1(X)$. Since we want to find an example of $U({\mathbb I};\
\mu;\ \e)$ satisfying Claim 1, we may assume, without loss of
generality, that $\mu(\{x\})=0,\ \forall x\in X$.
\\

\noindent {\bf Claim 3} {\it Let $U({\mathbb I}; \mu; \e)$ be given where
$\mu$ is a purely continuous measure on $X$. Let also $\overline V({\mathbb
I};\ \nu_1,...,\nu_n;\ \de)$ be a $\overline{\tau}$-neighborhood such that
$\mu$ and all $\nu_i$'s satisfy the condition
$$
\mu \ll \nu_1 \sim\cdots\sim\nu_n
$$
(i.e. $\mu$ is absolutely continuous with respect to all $\nu_i$'s which
are, in turn, pairwise equivalent). Then there exists $S\in \overline
V({\mathbb I};\ \nu_1,...,\nu_n;\ \de)$ but $S\notin U({\mathbb I};\ \mu;\
\e)$}.
\\

Let $f_1(x)=1$ and
$$
f_i(x)=\frac{d\nu_i}{d\nu_1(x)},\ \ \ i=2,...,n.
$$
Given $\de>0$, choose simple functions $g_i(x)$ where $g_1=1$,
$$
g_i=\sum_{j\in I_i}a_{ij}\chi_{E_{ij}},\ \ \ i=2,...,n,
$$
and such that
$$
\int_X \vert f_i-g_i\vert d\nu_1<\de/2,\ \ \ i=1,...,n.
$$
Define a new measure $\nu'_i$  on $X$ by
$d\nu'_i(x)=g_i(x)d\nu_i(x),\ i=1,2,...,n$.  Let ${\cal Q}_i$ be the
partition of $X$ defined by  $(E_{ij})_{j\in I_i}$. The
intersection of all ${\cal Q}_i$'s is a new partition ${\cal Q}=(F_{\bar
k})_{\bar k\in \La}$ of $X$ into Borel sets. Every $F_{\bar k}$ is of
the form $\bigcap_{i=1}^n E_{ij_i}$ where $j_i\in I_i$ and $\bar k$ is a
multiindex $(j_1,...,j_n)$ taken from a subset $\La \subset I_1\times
...\times I_n$. Let $S(\bar k): F_{\bar k} \to F_{\bar k}$ be a free
Borel map of $F_{\bar k}$ onto itself preserving $\nu_1$. Then for any
$F'\subset F_{\bar k}$ and $i=1,...,n$, we get
$$
\vert \nu'_i(F')-\nu'_i(S(\bar k)F')\vert =\vert
a_{ij_i}\nu_1(F')-a_{ij_i}\nu_1(S(\bar k)F')\vert =0,\ \ \bar k\in \La.
$$
It may be that some of the $F_{\bar k}$ have zero $\nu_1$-measure.
Thus, $S(\bar k)$ preserves the measures
$\nu_1'(=\nu_1),\nu_2',...,\nu_n'$, for all  $\bar k\in  \La$. Define
 a free Borel automorphism $S$ of $X$ by $Sx=S(\bar k)x$ when $x\in
F_{\bar k},\ \bar k\in \La$. Then $S\in \Aut$ and $Sx\ne x,\ x\in X$,
which means that $E(S,{\mathbb I})=X$ and $\mu(E(S,{\mathbb I}))=1$.
Moreover, $S$ preserves the measures $\nu_1'=\nu_1,\nu_2',...,\nu_n'$,
since for any Borel $F$, \begin{eqnarray*}
\nu'_i(SF) & = & \nu'_i(\bigcup_{\bar k\in \La}(SF\cap F_{\bar k}))\\
\\
& = &\sum_{\bar k \in \La}\nu'_i(S(F\cap F_{\bar k}))\\
\\
& = &\sum_{\bar k\in \La}\nu'_i(S(\bar k)(F\cap F_{\bar k}))\\
\\
& = & \sum_{\bar k\in \La}\nu'_i(F\cap F_{\bar k}))\\
\\
& = &\nu'_i(F).
\end{eqnarray*}
On the other hand , if $F\in {\cal B}(X)$, then

\begin{eqnarray*}
\vert \nu_i(F)-\nu_i(SF)\vert & = & \left|
\int_Ff_id\nu_1-\int_{SF}f_id\nu_1\right|\\
\\
& = & \left|
\int_F(f_i-g_i)d\nu_1+\int_Fg_id\nu_1-\int_{SF}(f_i-g_i)d\nu_1-\int_{SF}g_
id\nu_1\right|\\
\\
& \le & \int_F\vert f_i-g_i\vert d\nu_1+\int_{SF}\vert f_i-g_i\vert
d\nu_1+\left| \int_Fg_id\nu_1-\int_{SF}g_id\nu_1\right|\\
\\
& < & \de+\vert \nu'_i(F)-\nu'_i(SF)\vert\\
\\
& =& \de .
\end{eqnarray*}
Hence, we have shown that $S\in \overline V({\mathbb
I};\nu_1,...,\nu_n;\ \de)$ but $S\notin U({\mathbb I}; \mu; \e)$ if $\e
< 1$, and therefore Claim 3 is proved.
\\

\noindent
{\bf Claim 4} {\it Suppose $X_1\subset X$ is  a Borel set  such that
$\mu(X_1) > 0$ but $\nu_1(X_1) = \cdots =\nu_n(X_1) =0$. Then there exists
$S\in \Aut$ such that for any $\de > 0$, $S\in \overline V ({\mathbb I};\
\nu_1,...,\nu_n;\ \de)$ but $S\notin U({\mathbb I};\ \mu;\ \e)$ if $\e <
\mu(X_1)$.}
\\

To see that this claim holds, it is sufficient to take $S$ as a free
Borel automorphism on $X_1$ and put $Sx=x$ on $X - X_1$. Then
$\nu_i\circ S = \nu_i,\ i=1,...,n$, and therefore $S\in \overline V
({\mathbb I};\ \nu_1,...,\nu_n;\ \de)$ for any $\de > 0$. On the other
hand, $E(S,{\mathbb I}) = X_1$ and then $S$ cannot be in
$U({\mathbb I};\ \mu;\ \e)$ if $\e < \mu(X_1)$.
\\

$2^0$(c). Using Claims 3 and 4, we can prove the theorem in the general
case. To do this, it suffices to prove Claim 1. Let $U({\mathbb I};\ \mu;\
\e)$ be given where $\mu$ is a continuous measure on $X$. Let
$\nu_1,...,\nu_n$ be measures from ${\cal M}_1(X)$. Consider all possible
relations between $\nu_1,...,\nu_n$. For $\nu_1$ and $\nu_2$ there exists a
partition of $X$ into three Borel sets $A,B$, and $C$ such that $\nu_1$ and
$\nu_2$ are equivalent on $C$, $\nu_1$ is zero and $\nu_2$ is positive on
$B$, and $\nu_2$ is zero and $\nu_1$ is positive on $A$. Then $\nu_1$ is
supported on $A\cup C$ and $\nu_2$ is supported on $B\cup C$. Considering
the three measures $\nu_1,\nu_2$, and $\nu_3$, this fact can be applied to
each of the sets $A,B,C$. We thus obtain a new partition $ (X(1),\ X(2),\
X(3),\ X(1,2),\ X(1,3),\ X(2,3), \ X(1,2,3))$ of $X$ such that $\nu_i >0,\
\nu_j=0,\ j\neq i$ on $X(i)$, $\nu_i\sim\nu_k,\ \nu_j=0,\ j\neq i,k$ on
$X(i,k)$, and $\nu_1\sim\nu_2\sim\nu_3$ on $X(1,2,3)$. It is clear that a
similar statement holds for $\nu_1,...,\nu_n$. Namely, for every nonempty
subset $K\subset \{1,...,n\}$ there exists a subset $X(K)$ in $X$ such that
all the $\nu_i$'s are equivalent on $X(K)$ if $i\in K$, and $\nu_j(X(K))=0$
if $j\notin K$. Moreover, the sets $X(K)$ define a partition of $X$ as $K$
runs over all subsets in $\{1,...,n\}$.

Now consider the measure $\mu$ together with $\nu_1,...,\nu_n$. Without loss
of generality, we can assume that $\mu(X(K)) > 0$ for all $K$. Every $X(K)$
can be decomposed into two sets $X(K)'$ and $X(K)''$ (some of these sets
might be of zero $\mu$-measure) such that $\mu \ll \nu_i$ on $X(K)'$ and
$\mu $ is singular with respect to $\nu_i$ on $X(K)''$ where $i\in K$. In
fact, the latter condition holds for all $\nu_1,...,\nu_n$ by definition of
$X(K)$. This means that $\mu(X(K)'') > 0 $ whereas $\nu_i(X(K)'') = 0$,
$i=1,...,n$. Denote $X' = \bigcup_K X(K)',\ X'' = \bigcup_KX(K)''$. By Claim
4, we may define a free Borel automorphism $S$ on $X''$ which preserves each
$X(K)''$. To define $S$ on $X'$, we use the proof of Claim 3. Given
$\nu_1,...,\nu_n$ and $\de$, take $\de_1 < 2^{-n}\de$. The proof of Claim 3
allows us to find a Borel one-to-one map $S_K : X_K' \to X_K'$ such that
$S_Kx\neq x,\ x\in X(K)'$ and $\vert\nu_i(F) - \nu_I(S_KF)\vert < \de_1,\
F\subset X(K)'$. Then $Sx = S_Kx,\ x\in X(K)'$, defines a one-to-one Borel
map on $X'$. Therefore, $S\in \Aut$ and $S\in \overline V ({\mathbb I};\
\nu_1,...,\nu_n;\ \de)$ and since $S$ is free, $S\notin U({\mathbb I};\
\mu;\ \e)$ if $\e < 1$. The proof of Theorem \ref{T3.5} is
complete.\hfill{$\square$}

\begin{theorem}\label{T3.5.1} The topologies $p$ and $\tilde p$ are equivalent.
\end{theorem}

\noindent
{\it Proof}. Let $T\in \Aut$ and let $W(T)= W(T;F_1,...,F_n)$ be a
$p$-neighborhood of $T$. Then for every $S\in W(T)$, one has $SF_i =
TF_i,\ \forall i$. Let $f_i = \chi_{F_i}$. Then for any $0< \e <1$ we see
that $\widetilde W(T;f_1,...,f_n;\e) \subset W(T)$. Indeed, if
$$
\sup_{x\in X}|f_i(S^{-1}x) - f_i(T^{-1}x)| < \e,
$$
then $|\chi_{SF_i}(x) - \chi_{TF_i}(x)| < \e$ for all $x\in X$. This implies
that $SF_i = TF_i$.

To prove the converse statement, we take a $\tilde p$-neighborhood
$\widetilde W(T) = \widetilde W(T;f_1,...,f_n;\e)$. We need to show that
$\widetilde W(T)$ contains a $p$-neighborhood $W(T) = W(T; F_1,...,F_m)$.
Given $\e > 0$, find for each $f_i$ a Borel function $g_i(x)$ such that
$g_i(x) = \sum_{j\in I(i)} a_j(i)\chi_{E_j(i)(x)}$ and
$$
\sup_{x\in X} |f_i(x) - g_i(x)| < \e/2, \ i=1,...,n.
$$
Note that $|I(i)| < \infty$ since $f_i$ is bounded. Take the $p$-neighborhood
$W(T) = W(T; (E_j(i) : j\in I(i), i=1,...,n))$. If $S\in W(T)$, then
$$
\sup_{x\in X} |f_i(T^{-1}x) - f_i(S^{-1}x)|
$$
$$
\le \sup_{x\in X} |f_i(T^{-1}x) - g_i(T^{-1}x)|
+ \sup_{x\in X} |g_i(T^{-1}x) - g_i(S^{-1}x)| + \sup_{x\in X}
 |g_i(S^{-1}x) - f_i(S^{-1}x)|
$$
$$
< \e + \sup_{x\in X}|\sum_{j\in I(i)} a_j(i) (\chi_{TE_j(i)} -
\chi_{SE_j(i)}| = \e.
$$
Thus, $S\in \widetilde W(T)$. \hfill$\square$

\begin{proposition}\label{P3.6}  $(1)$ The topologies $\tau$ and $p$
are not comparable.\\ $(2)$ The topologies $\tau''$ and $\overline{p}$ are
not comparable.\\ (3) The topologies $\tau_0$ and $\tau''$ are not
comparable.
\end{proposition}

\noindent {\it Proof}. We prove only (3) here. We will show that for the
$\tau''$-neighborhood $U'' = U''({\mathbb I}; \delta_y; \e_0),\ e_0 <1$,
and any $\tau_0$-neighborhood $U_0 =U_0({\mathbb I};\nu_1,...,\nu_n; \e)$
there exists $T\in U_0$ such that $T\notin U''$. It suffices to take $T$
such that $\nu_i(E(T,{\mathbb I})) < \e$ for all $i$ and $Ty\neq y$.
Clearly such a $T$ can always be found. Then there is a Borel function
$f_0\in B(X)_1$ such that $|f_0(T^{-1}y) - f_0(y)| = 1$. Hence $T\notin
U''$. The fact that $\tau''$ cannot be stronger than $\tau_0$ is proved as
in Theorem \ref{T3.5}.

The proofs of (1) and (2) use the method similar to that applied in Claim 1
and in the proof of Theorem \ref{T3.5}. We leave the details to the reader.
\hfill$\square$

\begin{remark}\label{R3.7} {\rm Note that when the underlying space $X$ in
Definition \ref{D1.1} is a Polish space, one can consider various
modifications of the definition of $\tau, \tau', \tau''$ and $p,
\overline{p}$. In particular, we can show that replacing Borel functions
from $B(X)_1$ by continuous functions from $C(X)_1$ does not affect the
topology $\tau''$ (see (\ref{1.3})). To see this, we may use the following
statement in the proof of Proposition \ref{P3.4} which establishes the
equivalence of $\tau, \tilde\tau$, and $\overline{\tau}$:

{\it Let $g(x)\in B(X)_1$ and  $S\in \Aut$. Then for any $\de >0$ and
any $\mu \in {\cal M}_1(X )$, there exists a continuous function $f\in
C(X )_1$ such that}
$$
\left| \int_{X}(f(x)-g(x))d\mu \right| <\de,\ \ \ \ \ \ \ \ \ \left|
\int_{X}(f(S^{-1}x)-g(S^{-1}x))d\mu \right| <\de.
$$

In the case when $X$ is a Cantor set, observe that in the definition of
$\tau'$ (\ref{1.2}) it suffices to take the supremum over only the countable
family of clopen sets. This follows easily from regularity of Borel measures
on Cantor sets.}
\end{remark}


\sect{Bratteli diagram for a Borel automorphism}


\setcounter{equation}{0}

\noindent{\bf 4.1. Borel-Bratteli diagrams}. In this section, we show that
every aperiodic Borel automorphism of $\bs$ can be represented as a Borel
transformation acting on the space of infinite paths of a Bratteli diagram.
More precisely, we define a modification of the concept of Bratteli diagram
that is suitable for Borel automorphisms.

\begin{definition}\label{D4.1} A Borel-Bratteli diagram is an
infinite graph $B=(V,E)$ such that the vertex set $V$ and the edge
set $E$ can be partitioned into sets $V=\bigcup_{i\ge 0} V_i$ and
$E=\bigcup_{i\ge 1}E_i$ having  the following properties:

\noindent (i) $V_0=\{v_0\}$ is a single point, every $V_i$ and
$E_i$ are at most countable sets;

\noindent (ii) there exist a range map $r$ and a source map $s$
from $E$ to $V$ so that $r(E_i)\subset V_i$, $s(E_i)\subset
V_{i-1}$, $s^{-1}(v)\neq \emptyset $ for all $v\in V$,
and $r^{-1}(v)\neq\emptyset$ for all $v\in V \setminus V_0$.

\noindent (iii) for every $v\in V\setminus V_0$, the set $r^{-1}(v)$ is
finite.
\end{definition}

The pair $(V_i,E_i)$ will be called the $i$-th level. We write $e(v,v')$ to
denote an edge $e$ such that $s(e)=v,\ r(e)=v'$.

A finite or infinite sequence of edges, $(e_i : e_i\in E_i)$ such that
$s(e_i)=r(e_{i-1})$ is called a finite or infinite path, respectively. It
follows from the definition that every vertex $v\in V_i$ is connected to
$v_0$ by a finite path and the set of all such paths $E(v_0,v)$ is finite.
Given a Borel-Bratteli diagram $B = (V,E)$, we denote the set of infinite
paths by $Y_B$.

Let $B=(V,E,\geq)$ be a Borel-Bratteli diagram $(V,E)$ equipped with a
partial order $\geq$ defined on each $E_i,\ i=1,2,...,$ such that edges
$e,e'$ are comparable if and only if $r(e)=r(e')$; in other words, a linear
order $\geq$ is defined on each (finite) set $r^{-1}(v),\ v\in V\setminus
V_0$. For a Borel-Bratteli diagram $(V,E)$ equipped with such a partial
order $\geq$ on $E$, one can also define a partial lexicographic order on
the set $E_{k+1}\circ\cdots\circ E_l$ of all paths from $V_k$ to $V_l$:
$(e_{k+1},...,e_l) > (f_{k+1},...,f_l)$ if and only if for some $i$ with
$k+1\le i\le l$, $e_j=f_j$ for $i<j\le l$ and $e_i> f_i$. Then we see that
any two paths from $E(v_0,v)$, the (finite) set of all paths connecting
$v_0$ and $v$, are comparable with respect to the introduced lexicographic
order. We call a path $e= (e_1,e_2,..., e_i,...)$ maximal (minimal) if every
$e_i$ has a maximal (minimal) number amongst all elements from
$r^{-1}(e_i)$. Notice that there are unique minimal and maximal paths in
$E(v_0,v)$ for each $v\in V_i,\ i\ge 0$.

\begin{definition}\label{order} A Borel-Bratteli diagram $B=(V,E)$
together with a partial order $\geq$ on $E$ is called an ordered
Borel-Bratteli diagram $B=(V,E,\geq)$ if the space $Y_B$ has no cofinal
minimal and maximal paths. This means that $Y_B$ does not contain paths $e=
(e_1,..., e_i,...)$ such that for all sufficiently large $i$ the edges $e_i$
have maximal (minimal) number in the set $r^{-1}(e_i)$.
\end{definition}

It follows from the definitions that $Y_B$ is a 0-dimensional Polish space
in the natural topology defined by clopen sets.

To every Borel-Bratteli diagram $B = (V,E)$, we can associate a sequence of
infinite matrices. To do this, consider the set $E_n$ of all edges between
the levels $V_{n-1}$ and $V_n$. Let us enumerate the vertices of $V_n$ as
$(v_1(n),...,v_k(n),...), n\in \N$. Define the matrix $M_n =(a_{ik})_{i,k
=1}^\infty$ where $a_{ik} = |E(v_k(n-1),v_i(n))|$. We notice that only
finitely many entries of each row in $M_n$ are non-zero because of the
relation
$$
\sum_{k=1}^\infty a_{ik} = \vert r^{-1}(v_i^n)\vert < \infty.
$$
Moreover, in each column of $M_n$ there exists at least one non-zero entry.
Denote by ${\cal M}$ the set of infinite matrices having the above two
properties of their rows and columns. It is easy to see that ${\cal M}$ is
closed with respect to matrix multiplication.

Since the notion of ordered Bratteli diagram has been discussed in many
recent papers (see, e.g. [DH1], [DH2]), we recall only the principal
definitions. We refer to [GPS, HPS] for more detailed expositions of this
material. We will also use in our proofs the notions of {\it telescoping and
splitting} of a Bratteli diagram defined in [GPS]. By definition, the
telescoping $B' = (V',E', \geq)$ of $B =(V,E, \geq)$ with respect to a
sequence $0= m_0 < m_1 < m_2 < \cdots$ is obtained if we set $V'_n =
V_{m_n}$ and $E'_n= E_{m_{n-1}+1}\circ\cdots\circ E_{m_n}$. The set $E'_n$
has the induced lexicographic order defined above. The operation converse to
telescoping a Bratteli diagram is splitting. This means that between to
consecutive levels, say $V_{n-1}$ and $V_n$, we add a new level $V'$ which
is a disjoint union of finite sets $V'(v)$ such that the number of vertices
in $V'(v)$ equals the number of edges in $r^{-1}(v),\ v\in V_n$. It is easy
to see how to introduce an order on the edge set of the new diagram so that
by telescoping one gets the original ordered diagram back.

For each ordered Borel-Bratteli diagram $B=(V,E,\geq)$, define a Borel
transformation $\varphi$ (also called the Vershik automorphism) acting on
$Y_B$ as follows. Given $y=(e_1, e_2,...)\in Y_B$, let $k$ be the smallest
number such that $e_k$ is not a maximal edge. Let $f_k$ be the successor of
$e_k$ in $r^{-1}(r(e_k))$. Then we define $\varphi(x)=
(f_1,...,f_{k-1},f_k,e_{k+1},...)$ where $(f_1,..., f_{k-1})$ is the minimal
path in $E(v_0, r(f_{k-1}))$. Obviously, $\varphi$ is a one-to-one mapping
of $Y_B$ onto itself. Moreover, $\varphi$ is a homeomorphism of $Y_B$.

It is sometimes convenient to use another description of infinite paths in
an ordered Borel-Bratteli diagram $B= (V,E,\geq)$. Take $v\in V_n$ and
consider the finite set $E(v_0,v)$. The lexicographic order on $E(v_0,v)$
gives us an enumeration of its elements from 0 to $h(n,v)-1$ where 0 is
assigned to the minimal path and $h(n,v)-1$ is assigned to the maximal path
in $E(v_0,v)$. Note that $h(1,v) = |r^{-1}(v)|,\ v\in V_1$ and by induction
\begin{equation}\label{4.1}
h(n,v) = \sum_{w\in s(r^{-1}(v))}|E(w,v)|h(n-1,w),\ \ \ v\in V_n.
\end{equation}

Let $y=(e_1,e_2,...)$ be an infinite path from $Y_B$. Consider a sequence
$(P_n)$ of enlarging finite paths defined by $y$: $P_n = (e_1,...,e_n)\in
E(v_0, r(e_n))$, $n\in \N$. Then every $P_n$ can be identified with a pair
$(i_n,v_n)$ where $v_n = r(e_n)$ and $i_n \in [0, h(n,v_n)-1]$ is the number
assigned to $P_n$ in $E(v_0,v_n)$. Thus, every $y = (e_n)\in Y_B$ can be
uniquely represented as an infinite sequence $(i_n,v_n)$ with $v_n =r(e_n)$
and $0\leq i_n \leq h(n,v_n)-1$. We also observe that $i_n \to \infty,\ n\to
\infty, $ since the Borel-Bratteli diagram $B$ has no infinite cofinal
minimal paths and $(h(n,v_n) - i_n) \to \infty,\ n\to\infty,$ since there is
no infinite cofinal maximal paths in $B$.

Thus, given an ordered Borel-Bratteli diagram $B=(V,E,\geq)$, we have
defined a dynamical system $(Y_B, \varphi)$. Our goal is to show that every
Borel automorphism can be realized as a Vershik transformation acting on the
space of infinite paths of an ordered Borel-Bratteli diagram.
\\

\noindent {\bf 4.2. Construction of a Borel-Bratteli diagram by a Borel
automorphism and a vanishing sequence of markers}. Let $T$ be an aperiodic
automorphism of $(X,\B)$. Take a vanishing sequence of markers $(A_n)$ with
$X =A_0 \supset A_1\supset A_2\supset\cdots$ (see Definition \ref{markers}).
We give a construction of an ordered Borel-Bratteli diagram coming from
$(A_n)$ and $T$.

Since $A_1$ is a recurrent complete $T$-section of $X$, there exists a
partition $\xi_1= \bigcup_{v\in V_1}\xi_1(v)$ of $X$ formed by at most
countable collection of disjoint $T$-towers $\xi_1(v) = (A_1(v),
T(A_1(v)),\dots , T^{h(1,v)-1}(A_1(v)))$ where $\bigcup_{v\in V_1}A_1(v) =
A_1$ and $V_1$ is a subset of $\N$. Here $h(1,v)$ is the height of $T$-tower
$\xi_1(v)$.

Because $A_2$ is a subset of $A_1$, we can assume (refining, if necessary,
the partition $\xi_1$) that $A_2$ is a union of some sets $A_1(v)$. Then we
again define $\xi_2$ as a disjoint collection of $T$-towers $\xi_2(v) =
(A_2(v),T(A_2(v)),\dots ,T^{h(2,v)-1}(A_2(v))),\ v\in V_2\subset \N,$ with
$\bigcup_{v\in V_2}A_2(v) = A_2$. We apply this construction for every $A_n$
and find the corresponding partition $\xi_n$ consisting of $T$-towers
$\xi_n(v)$ of finite height $h(n,v)$, $v\in V_n\subset \N$. Note that at
each step, $A_{n+1}$ is a $\xi_n$-set and hence $\xi_{n+1}$ refines $\xi_n$.
Moreover, for any partition $\xi$ from the sequence $(\xi_n)$ and any Borel
set $D$ in $X$ we can refine $\xi$ such that $D$ becomes a $\xi$-set. This
means that we can assume that the collection of atoms of $(\xi_n)$ separates
points in $X$.

Now define an ordered Borel-Bratteli diagram $B$ using the sets $A_n$ and
generated by them partitions $\xi_n$ of $X$, $n\in \N$. Let $V_0 =\{v_0\}$
be a singleton (relating to $A_0 =X$). The set $V_1$ gives vertices at the
first level in $B$. To define $E_1$, we take $v\in V_1$ and draw $h(1,v)$
edges connecting $v_0 $ and $v$. Enumerate these edges from 1 to $h(1,v)$ in
an arbitrary order. Set $s(e) = v_0, r(e)= v$ for $e$ connecting $v_0$ and
$v$. Thus, the set $r^{-1}(v)$ becomes linearly ordered for every $v\in
V_1$.

To define the diagram $B$ for the next level, we take $V_2$, obtained from
the partition $\xi_2$, as the set of vertices. Fix a vertex $v \in V_2$ and
consider the $T$-tower $\xi_2(v)$. It can be easily seen that $\xi_2(v)$
intersects a finite number of $T$-towers from the partition $\xi_1$, say
$\xi_1(v_1), \xi_1(v_2),..., \xi_1(v_s),\ v_1,...,v_s\in V_1$. Notice that
these towers are not necessarily different and that some of the $v_i$'s may
be met several times. We see that $h(2,v) = h(1,v_1) +\cdots + h(1,v_s)$.
Take the vertices $v_1,...,v_s$ from $V_1$ and draw the edges connecting
each of them with $v\in V_2$. We get the sets $E(v_i,v),\ i=1,...,s$. The
number of edges in $E(v_i,v)$ equals the multiplicity of the vertex $v_i$
in the set $v_1,...,v_s$. Define $s(e)= v_i, r(e) = v$ where $e\in
E(v_i,v)$. To introduce a linear order on $r^{-1}(v),\ v\in V_2$ we
consider again the $T$-towers $\xi_1(v_1), \xi_1(v_2),..., \xi_1(v_s)$ and
enumerate them from 1 to $s$ according to the natural order in which they
appear in $\xi_2(v)$ when we go along $\xi_2(v)$ from the base to the top.
Enumerate the corresponding edges $e\in r^{-1}(v)$ in the same order from 1
to $s$. This procedure is applied to every vertex from $V_2$ to define the
entire set $E_2$ together with a partial order on $E_2$.

Repeating this method for every $n$, we construct the $n$-th level $(V_n,
E_n)$ and establish a partial order on $E_n$. We see that conditions (i) -
(iii) of Definition \ref{D4.1} hold for the infinite graph $B = (V,E)$ where
$V=\bigcup_{n\ge 0} V_n$ and $E=\bigcup_{n\ge 1} E_n$ are defined as above.
The partial order which we have described on $E$ determines an ordered
Borel-Bratteli diagram $B = (V,E,\geq)$ according to Definition \ref{order}.
Indeed, it is easy to see that $B$ has no cofinal maximal and minimal paths.
It follows from the fact that $\bigcap_n T^k(A_n)= \emptyset$ for any $k\in
\Z$. We also observe that every infinite path $y\in Y_B$ is completely
determined by the infinite sequence $\{(i_n,v_n)\}_n,\ v\in V_n, 0\leq i\leq
h(n,v_n)-1$ such that $T^{i_{n+1}}(A_{n+1}(v_{n+1})) \subset
T^{i_{n}}(A_{n}(v_{n})),\ n\in \N$. The height $h(n,v)$ can be found by
(\ref{4.1}).

Notice that if one takes a subsequence $(A_{n_m})$ in $(A_n)$ and constructs
a new ordered Borel-Bratteli diagram $B'$ by $(A_{m_n})$ and $T$, then $B'$
turns out to be a telescoping of $B$. If a set $A_n(v)$ is partitioned into
a finite number of uncountable Borel sets $A^v_n(w)$ and respectively
$\xi_n(v)$ is cut into a finite number of new $T$-towers $\xi_n^v(w)$ with
base $A_n^v(w)$, then the above construction gives an ordered Borel-Bratteli
diagram $B''$ which is a splitting of $B$.
\\

The next theorem shows that any aperiodic Borel automorphism $T$ can be
realized as a Vershik transformation acting on the space of infinite paths
of an ordered Borel-Bratteli diagram.

\begin{theorem}\label{T4.2} Let $T$ be an aperiodic Borel automorphism acting
on a standard Borel space $(X, \B)$. Then there exists an ordered
Borel-Bratteli diagram $B=(V,E,\geq)$ and a Vershik automorphism $\varphi :
Y_B \to Y_B$ such that $(X, T)$ is isomorphic to $(Y_B,\varphi)$.
\end{theorem}

\noindent {\it Proof}. Let $(A_n)$ be a vanishing sequence of markers for
$T$ and let $\xi_n = (\xi_n(v): v\in V_n)$ be a collection of disjoint
$T$-towers where $\xi_n(v) = (A_n(v), T(A_n(v)),...,
T^{h(n,v)-1}(A_n(v))),\ n\in \N$. As mentioned above, we can assume that
the atoms of $(\xi_n,\ n\in \N)$ generate the Borel structure on $X$. By
changing-of-topology results (see Remark \ref{change}), we may choose a
topology $\omega$ on $X$ such that: (i) $X$ is a Polish 0-dimensional
space, (ii) $\B(\omega)= \B$ where $\B(\omega)$ is the $\sigma$-algebra
generated by $\omega$-open sets, (iii) all sets $T^j(A_n(v)),\ v\in V_n,
j=0,1,...,h(n,v)-1, n\geq 1,$ are clopen in $\omega$, (iv) $T$ is a
homeomorphism of $(X,\omega)$.

Next, we observe that for fixed $n\in\N$ and $\e >0$ one can cut each
$T$-tower $\xi_n(v),\ v\in V_n,$ into disjoint clopen towers of the same
height such that the diameter of every element of the new towers is less
than $\e$. Therefore, without loss of generality, we may assume that
\begin{equation}\label{4.2}
\sup_{0\leq j< h(n,v), v\in V_n}[{\rm diam}\ T^j(A_n(v))]\to 0,\ \ n\to
\infty.
\end{equation}

Now applying the above construction to $(A_n)$ and $T$ (see subsection {\bf
4.2}), we can find an ordered Borel-Bratteli diagram $B=(V,E,\geq)$
together with a Vershik transformation acting on the space of infinite
paths $Y_B$. Define a map $F: X\to Y_B$. Given $x\in X$, choose the unique
sequence $\{(i_n,v_n)\}_n,\ 0\leq i\leq h(n,v_n)-1, v_n\in V_n$, such that
\begin{equation}\label{4.3}
T^{i_{n+1}}(A_{n+1}(v_{n+1})) \subset T^{i_{n}}(A_{n}(v_{n}))
\end{equation}
and $\{x\} = \bigcap_n T^{i_n}(A_n(v_n))$. As noticed in section {\bf
4.1}, such a sequence defines a unique infinite path $y\in Y_B$. We set
$F(x) = y$. It is clear that $F$ is a continuous injection of $X$ into
$Y_B$. Moreover, it easily follows from the construction of $B$ and
definition of Vershik transformation $\varphi$, acting on $Y_B$, that
$F(Tx) = \varphi(Fx),\ x\in X$. To finish the proof, we need to show only
that $F(X) = Y_B$. Take a path $y\in Y_B$. Then $y$ defines an infinite
sequence $\{(i_n,v_n)\}_n,\ v\in V_n, 0\leq i\leq h(n,v_n)-1$ as in {\bf
4.1}. By (\ref{4.2}) and (\ref{4.3}) we get a single point $x$ such that
$F(x) = y$. \hfill$\square$

\begin{remark}\label{R4.3} {\rm  The idea of the proof of Theorem
\ref{T4.2} was shown to us by B.~Miller [M]. He has pointed out that if
(\ref{4.2}) does not hold, then $F$ may be a map onto a proper subset of
$Y_B$. However, one can show that in any case $F(X)$ contains a dense
$G_\delta$-set. Indeed, define
$$
D= \bigcap_{n=1}^\infty\bigcap_{m=1}^\infty\bigcup_{k\geq
n}^\infty\bigcup_{v\in I(m,k)} A_k(v)
$$
where $I(m,k) = \{v\in V_k : {\rm diam}\ (A_k(v)) < 1/m\}$. Then $F(D)$ is a
dense $G_\delta$-subset of $Y_B$.}
\end{remark}

We observe that given aperiodic $T\in \Aut$ one can find a vanishing
sequence of markers such that the corresponding ordered Borel-Bratteli
diagram has a finite number of vertices at each level. To do this, we
follow a suggestion of B.~Miller [M] using maximal sets.

By definition, an uncountable Borel set $A$ is called {\it $k$-maximal} for
an aperiodic automorphism $T\in \Aut$ if $A\cap T^iA = \emptyset,\
i=1,...,k-1$, and $A$ cannot be extended to a larger set having this
property. It is easy to show that $A$ is $k$-maximal if and only if
\begin{equation}\label{4.4}
X =\bigcup_{|i|< k}T^iA\ \ {\rm and}\ \ A\cap T^iA =\emptyset,\ i=1,...,k-1.
\end{equation}
Hence, a maximal set $A$ is a complete section for $T$ such that every point
from $A$ is recurrent. The existence of maximal sets can be easily deduced,
for instance, from the tower construction used in Section 2. Indeed, let
$\xi = (\xi(v) : v\in V)$, where $\xi(v) = (B(v), TB(v),..., T^{h(v)-1}B)$,
be a partition of $X$ into $T$-towers. Define
\begin{equation}\label{4.4'}
A= \bigcup_{v\in V}\bigcup_{i=0}^{h(v)k^{-1} -1}T^{ik}B(v).
\end{equation}
It follows from (\ref{4.4}) that $A$ is a $k$-maximal set.

\begin{proposition}\label{P4.4} Given aperiodic $T\in \Aut$, there exists a
vanishing sequence of markers $(A_i)$ such that the corresponding ordered
Borel-Bratteli diagram has a finite number of vertices at each level.
\end{proposition}

\noindent{\it Proof}. Let $A_0 = X$ and let $A_1$ be a $k$-maximal set for
$T$. Since every point $x\in A_1$ returns to $A_1$, one can define the
induced aperiodic Borel automorphisms $T_{A_1}$ acting on $A_1$:
$$
T_{A_1}x = T^{m(x)}x,\ \ {\rm where}\ \ m(x) = \min\{m>0 : T^m x\in A_1\}.
$$
Then $T_{A_1}$ is again an aperiodic automorphism of $A_1$ and we can find
an $k$-maximal set $A_2\subset A_1$ for $T_{A_1}$. Let now $A_{i+1}\subset
A_i$ be a $k$-maximal set for $T_{A_i},\ i\in \N$. It follows from
(\ref{4.4}) that the construction used in subsection {\bf 4.2} gives a
finite number of $T$-towers over $A_1$ which cover $X$. By the same
reasoning, we see that $A_1$ is covered by a finite number of disjoint
$T_{A_1}$-towers constructed over $A_2$. Therefore, $X$ is covered by a
finite number of disjoint $T$-towers constructed over the base $A_2$. It is
easy to see that we will have this property at every stage of the
construction.

The sequence $(A_i)_{i\in \N}$ of decreasing $k$-maximal sets which we have
defined, may have a non-empty intersection, $A_\infty = \bigcap_{i\in \N}
A_i$. Obviously, $A_\infty$ is a wandering set with respect to $T$. We set
$A'_i = (A_i \setminus A_\infty)\cup (\bigcup_{|j|>i} T^jA_\infty)$. Then
$A'_i\supset A'_{i+1}$, $\bigcap_iA'_i =\emptyset$, and therefore $(A'_i)$
is a desired vanishing sequence of markers. \hfill$\square$
\\

The concept of Borel-Bratteli diagram can be used to obtain another proof of
the following result (see [N, Theorem 8.9]).

\begin{corollary}\label{C4.6} If $T$ is an aperiodic homeomorphism of a
Polish space $X$, then there exists a compact metric space $Y$ and a
homeomorphism $S$ of $Y$ such that $T$ is homeomorphic to the restriction of
$S$ to an $S$-invariant dense $G_\delta$-subset of $Y$.
\end{corollary}

\noindent{\it Proof}. We will use the method of proof of Theorem \ref{T4.2}
and Proposition \ref{P4.4} to construct a vanishing sequence of markers
$(A_n)$ such that the corresponding ordered Borel-Bratteli diagram $B$ has
a finite number of vertices at each level and such that every sequence
satisfying (\ref{4.3}) has a non-empty intersection.

We start with a vanishing sequence of markers $(D_n)$ satisfying
(\ref{4.2}). Let $\xi_1$ be a partition of $X$ into $T$-towers constructed
over $D_1$. Define the $k$-maximal set $A_1$ as in (\ref{4.4'}). Clearly,
$A_1\supset D_2$. Take the induced automorphism $T_{A_1}$ and construct the
partition $\xi_2$ of $A_1$ into $T_{A_1}$-towers over $D_2$. Let $A_2$ be a
$k$-maximal subset of $A_1$ defined again as in (\ref{4.4'}). Continuing
this process we define a decreasing sequence of Borel sets $(A_n)$. Let
$(\eta_n)$ be the partitions of $X$ into $T$-towers constructed over the
$A_n$'s. By Proposition \ref{P4.4}, the corresponding ordered
Borel-Bratteli diagram $B$ has a finite number of vertices at each level.
Moreover, we can assume that the atoms of the partitions $\eta_n, n\in \N$,
separate points of $X$ (see subsection {\bf 4.2}). If
$\{T^{i_{n}}(A_{n}(v_{n}))\}$ is a sequence of atoms satisfying
(\ref{4.3}), then the intersection of these atoms contains at most one
point. But, in fact, this intersection is non-empty because it contains the
intersection of a similar sequence of atoms of partitions defined by
$(D_n)$.

As in the proof of Theorem \ref{T4.2}, take the topology $\omega$ satisfying
conditions (i) - (iv). Clearly, the ordered Borel-Bratteli diagram $B$ has
maximal and minimal paths (not necessarily unique). To define a
homeomorphism $S$ of a compact space $Y$, we use Forrest's construction of
path-sequence dynamical system generated by a Bratteli diagram [For]. By
definition, the set $Y = Y_1\cup Y_2$ where $Y_1$ is the set of all infinite
paths from $Y_B$ different from orbits of cofinal maximal and minimal paths
and $Y_2$ is the set of all equivalent pairs in the sense of Forrest $(x,y)$
with $x$ maximal and $y$ minimal paths in $Y_B$. Let $S$ be the
homeomorphism defined in [For] on $Y$. Then $(Y,S)$ is a Cantor dynamical
system and the set $Y_1$ is an $S$-invariant dense $G_\delta$-subset of $Y$.
There is a homeomorphism between the action of $S$ on $Y_1$ and the action
of $T$ on $X$. \hfill$\square$

\begin{remark} {\rm Suppose that a Borel-Bratteli diagram $B= (V,E,\geq)$ has
a finite number of vertices at each level, i.e. $|V_n| < \infty,\ n\in \N$.
Then either (i) $\limsup_n |V_n| = K <\infty$ or (ii) $\limsup_n |V_n| =
\infty$. If (i) holds, then there exists a subsequence $(n(k))$ such that
$|V_{n(k)}| = K$ for all $k$. By telescoping, we can produce a new
Borel-Bratteli diagram $B'$ isomorphic to $B$ such that the number of
vertices at each level of $B'$ is exactly $K$. If (ii) holds, then there
exists a subsequence $(n(k))$ such that $|V_{n(k)}| < |V_{n(k+1)}|$ for all
$k$. This means that for the Borel-Bratteli diagram obtained by telescoping
with respect to $(n(k))$, the number of vertices at each level is a
strictly monotonically increasing sequence.}
\end{remark}

Another application of Proposition \ref{P4.4} consists of description of
the closures of odometers. We recall that in the case of Cantor dynamics
the closure of odometers in $D$ coincides with the set of moving
homeomorphisms (by definition, $S$ is moving if for any proper clopen set
$E$ the sets $E\setminus SE$ and $SE\setminus E$ are not empty). In Borel
dynamics the notion of moving automorphisms has no sense.

\begin{theorem} \label{od} $(1)$\ $\overline{{\cal O}d}^\tau = \ap$ and
$\overline{{\cal O}d}^{\tau_0} = \ap \mod(Ctbl)$.\\
$(2)$\ $\overline{{\cal O}d}^p \subset \inc$.\\
$(3)$\  $\overline{{\cal S}m}^D \supset \overline{{\cal O}d}^D \supset \ap$
 assuming that $(X,d)$ is a compact metric space.\\
\end{theorem}
{\it Proof}. (1) Given $T\in \ap$ and $\e> 0$, choose a natural number $n >
1/\e$. By Proposition \ref{P4.4}, we can find a partition $\xi$ of $X$ into
a finite number of disjoint  $T$-towers $\xi(v)= (A(v), TA(v),...,
T^{h(v)-1}A(v)),\ v\in V,\ |V| < \infty$, such that $\min (h(v) : v\in V)
\geq 2n$.  We call the $\xi$-set
$$
L(i) = \bigcup_{v\in V} T^{h(v)-i-1} A(v),\ \  i=0,1,...,n-1,
$$
the $i$-th level in the partition $\xi$. Let $\mu$ be a Borel probability
measure on $X$. Then there is a pair $(L(i_0), L(i_0+1))$ such that
$\mu(L(i_0)\cup L(i_0 +1)) < \e, i_0 = 0,1,...,n-1$. Take now the set
$L(i_0)$ and construct a new partition $\xi'$ of $X$ into disjoint
$T$-towers $\xi'(j),\ j = 1,...,k$, with the base $L(i_0)$ and the top
$L(i_0+1)$. Define the automorphism $S$ of $X$ to coincide with $T$
everywhere except on the top level $L(i_0+1)$. We define $S$ on $L(i_0+1)$
as a Borel one-to-one map from the top of the tower $\xi'(j)$ onto the base
of $\xi'(j+1),\ j-1,...,k-1,$ and the top of $\xi'(k)$ is sent by $S$ onto
the base of $\xi'(1)$. In such a way, the space $X$ is represented as an
$S$-tower. It is easily seen that the definition of $S$ can be refined to
produce an odometer $S_1$ which agrees with $S$ everywhere except the top
of the $S$-tower. By construction, $\mu(E(S_1, T)) < \e$. The fact that the
set ${\cal O}d$ is dense in $\ap$ with respect to $\tau$ follows now from
the latter inequality and Remark \ref{R1.4} (4).

To prove the other statement of (1), we note first that by Theorem
\ref{T2.6}
$$
\overline{{\cal O}d}^{\tau_0} \subset \ap \mod(Ctbl)
$$
since ${\cal O}d \subset \ap$ and
$\overline{\ap}^{\tau_0} = \ap \mod(Ctbl)$. On the other hand, the topology
$\tau$ is stronger than $\tau_0$ and therefore
$$
\overline{{\cal O}d}^{\tau_0} \supset \overline{{\cal O}d}^{\tau}
= \ap
$$
Hence $\ap \mod(Ctbl)\subset \overline{{\cal O}d}^{\tau_0}$
and we are done.

(2) This follows from the fact that $\inc$ is closed in $p$ (Theorem
\ref{inc}).

(3) We note  that once we have a finite partition of $X$ into $T$-towers as
in Proposition \ref{P4.4}, then every tower can be additionally cut into
finitely (or countably) many subtowers such that the diameter of every atom
is sufficiently small. To construct either an odometer or a smooth
automorphism, we can use the same method as in the proof of (1).
\hfill$\square$

\noindent{\bf 4.3. Special Borel-Bratteli diagrams}. In this subsection we
first define the notion of special Borel-Bratteli diagrams and then
indicate a class of automorphisms which are completely described by these
diagrams.

By definition, an ordered Borel-Bratteli diagram $B=(V,E,\geq)$ is called
{\it special} if it satisfies the following conditions:

(i) $V=\bigcup_{n\geq 0} V_n$ and for $n\geq 1$, $V_n= \bigcup_{j\geq
n}V_{nj}$, where $V_{nj}\cap V_{nj'} =\emptyset,\ j\neq j'$, and $2 \leq
|V_{nj}| < \infty$. The set $V_{nn}$ is a union of two disjoint sets
$V_{nn}(0)$ and $V_{nn}(1)$ with $|V_{nn}(0)| \geq 2$.

(ii) $E= \bigcup_{n\geq 1}E_n$ where each $E_n$ is a union of disjoint
finite subsets $E_{nj},\ j \geq n$, such that for $j > n$ and $e\in E_{nj}$,
one has $s(e) \in V_{n-1,j},\ r(e) \in V_{nj}$ and $s(E_{nj}) = V_{n-1,j}$.
The set $E_{nn}$ is a union of two disjoint subsets $E_{nn}(0)$ and
$E_{nn}(1)$ such that $s(e) \in V_{n-1,n-1},\ r(e) \in V_{nn}(0)$ if $e\in
E_{nn}(0)$ and $s(e) \in V_{n-1,n},\ r(e) \in V_{nn}(1)$ if $e\in
E_{nn}(1)$. Moreover $|r^{-1}(v)| =1$ if $v$ is either in $V_{nj},\ j >n$,
or $v\in V_{nn}(1)$. If $v\in V_{nn}(0)$, then $|r^{-1}(v)| \geq 4$. The
edges $e_1 <\cdots < e_m$ from $|r^{-1}(v)|,\ v\in V_{nn}(0)$ are ordered
such that $s(e_1) \in V_{n-1,n},\ s(e_i) \in V_{n-1,n-1},\ i=2,...,m-1,$ and
$s(e_m)$ is either in $V_{n-1,n}(1)$ or in $V_{n-1,j},\ j>n-1$. By
definition, $s(E_{nj})= V_{n-1,j}$.

Figure 1 is an example of a special Borel-Bratteli diagram which illustrates
the above definition (we do not indicate a partial order on edges of this
diagram). It follows from the definition (see also Figure 1) that the set of
infinite paths does not have cofinal minimal and maximal paths.
\\

\unitlength=1.00mm \special{em:linewidth 0.4pt} \linethickness{0.4pt}
\begin{picture}(137.33,145.67)
\put(70.00,140.00){\circle*{1.50}} \put(2.00,116.00){\circle*{1.50}}
\put(12.00,116.00){\circle*{1.50}} \put(27.00,116.00){\circle*{1.50}}
\put(37.00,116.00){\circle*{1.50}} \put(67.00,116.00){\circle*{1.50}}
\put(77.00,116.00){\circle*{1.50}} \put(107.00,116.00){\circle*{1.50}}
\put(117.00,116.00){\circle*{1.50}} \put(127.00,116.00){\circle*{0.50}}
\put(132.00,116.00){\circle*{0.50}} \put(137.00,116.00){\circle*{0.50}} \
\put(2.00,80.00){\circle*{1.50}} \put(12.00,80.00){\circle*{1.50}}
\put(27.00,80.00){\circle*{1.50}} \put(37.00,80.00){\circle*{1.50}}
\put(67.00,80.00){\circle*{1.50}} \put(77.00,80.00){\circle*{1.50}}
\put(107.00,80.00){\circle*{1.50}} \put(117.00,80.00){\circle*{1.50}}
\put(127.00,80.00){\circle*{0.50}} \put(132.00,80.00){\circle*{0.50}}
\put(137.00,80.00){\circle*{0.50}} \ \ \ \put(2.00,45.00){\circle*{1.50}}
\put(12.00,45.00){\circle*{1.50}} \put(27.00,45.00){\circle*{1.50}}
\put(37.00,45.00){\circle*{1.50}} \put(67.00,45.00){\circle*{1.50}}
\put(77.00,45.00){\circle*{1.50}} \put(107.00,45.00){\circle*{1.50}}
\put(117.00,45.00){\circle*{1.50}} \put(127.00,45.00){\circle*{0.50}}
\put(132.00,45.00){\circle*{0.50}} \put(137.00,45.00){\circle*{0.50}} \
\put(62.00,25.00){\circle*{0.50}} \put(72.00,25.00){\circle*{0.50}}
\put(82.00,25.00){\circle*{0.50}} \ \
\bezier{308}(1.67,115.67)(32.33,140.33)(70.33,140.00)
\bezier{172}(37.00,115.67)(59.33,123.33)(70.00,140.00)
\put(12.00,115.67){\line(5,2){58.00}} \put(27.33,116.00){\line(2,1){42.00}}
\put(66.67,116.00){\line(1,6){3.67}} \put(77.00,116.33){\line(-1,4){5.67}}
\put(106.67,116.00){\line(-3,2){34.33}}
\put(117.33,116.00){\line(-2,1){43.67}}
\bezier{144}(1.33,79.67)(0.00,98.67)(1.00,116.00)
\bezier{300}(1.67,80.00)(40.00,97.00)(67.00,115.33)
\bezier{152}(12.33,79.67)(2.00,95.00)(1.67,114.67)
\put(2.00,80.00){\line(2,3){24.33}}
\bezier{80}(0.00,116.00)(9.67,121.67)(16.33,115.67)
\bezier{72}(0.33,116.00)(7.00,111.33)(16.00,115.67)
\bezier{88}(23.33,115.67)(33.33,121.67)(42.00,116.00)
\bezier{88}(23.00,116.00)(31.67,111.00)(42.33,115.67)
\bezier{92}(62.00,115.67)(71.67,121.33)(82.33,116.00)
\bezier{88}(61.67,115.67)(72.33,112.33)(82.67,115.67)
\bezier{88}(101.67,115.67)(111.00,121.33)(121.00,116.00)
\bezier{96}(101.67,115.67)(113.33,109.33)(122.33,115.67)
\put(70.00,145.67){\makebox(0,0)[cc]{$V_0$}}
\put(2.00,123.67){\makebox(0,0)[cc]{$V_{11}(0)$}}
\put(45.67,115.00){\makebox(0,0)[cc]{$V_{11}(1)$}}
\put(87.33,115.33){\makebox(0,0)[cc]{$V_{12}$}}
\put(124.00,111.00){\makebox(0,0)[cc]{$V_{13}$}}
\put(67.33,79.33){\line(6,5){38.33}} \put(77.00,80.33){\line(6,5){38.00}}
\put(26.33,79.67){\line(6,5){40.00}} \put(36.33,80.00){\line(6,5){37.67}}
\bezier{96}(102.67,79.33)(112.67,87.00)(121.67,79.67)
\bezier{88}(102.67,78.67)(113.33,74.67)(122.33,79.67)
\bezier{92}(62.33,79.33)(74.33,86.00)(81.00,79.67)
\bezier{80}(63.00,79.33)(72.00,74.67)(81.00,79.00)
\bezier{88}(23.67,79.67)(33.33,86.33)(41.33,79.33)
\bezier{76}(23.67,79.67)(31.67,75.67)(41.33,79.00)
\bezier{72}(0.33,80.00)(6.67,84.67)(15.67,79.67)
\bezier{68}(0.00,79.33)(7.00,76.00)(15.67,79.67)
\bezier{76}(0.00,44.67)(7.00,49.67)(15.67,44.67)
\bezier{72}(0.00,44.33)(7.00,39.67)(15.33,44.00)
\bezier{80}(23.00,45.00)(32.67,49.67)(41.00,45.00)
\bezier{80}(23.00,44.00)(31.00,40.67)(41.33,44.33)
\bezier{92}(62.33,45.00)(73.00,51.00)(82.33,45.00)
\bezier{84}(62.33,44.33)(71.67,41.00)(82.67,44.33)
\bezier{88}(103.33,44.67)(111.33,51.00)(121.00,45.00)
\bezier{84}(103.67,44.67)(113.67,38.67)(121.00,44.67)
\put(107.00,80.00){\line(5,6){21.33}} \put(116.67,80.67){\line(3,4){19.33}}
\put(107.00,44.33){\line(4,5){24.33}} \put(117.67,44.33){\line(2,3){19.67}}
\put(66.33,45.00){\line(6,5){40.00}} \put(77.00,44.67){\line(6,5){39.00}}
\put(26.00,44.67){\line(6,5){39.67}} \put(37.00,44.67){\line(6,5){37.67}}
\put(2.00,80.33){\line(1,1){33.33}} \put(12.00,79.67){\line(0,1){36.33}}
\put(12.00,79.67){\line(5,3){54.00}} \put(12.33,45.00){\line(-1,3){11.33}}
\put(11.67,44.67){\line(2,5){14.00}} \put(11.33,44.00){\line(3,2){53.67}}
\put(10.67,73.67){\makebox(0,0)[cc]{$V_22(0)$}}
\put(44.67,76.00){\makebox(0,0)[cc]{$V_{22}(1)$}}
\put(85.33,77.67){\makebox(0,0)[cc]{$V_{23}$}}
\put(120.33,75.33){\makebox(0,0)[cc]{$V_{24}$}}
\put(11.33,37.33){\makebox(0,0)[cc]{$V_{33}(0)$}}
\put(43.33,38.00){\makebox(0,0)[cc]{$V_{33}(1)$}}
\put(86.00,40.67){\makebox(0,0)[cc]{$V_{34}$}}
\put(120.67,37.67){\makebox(0,0)[cc]{$V_{35}$}}
\put(2.00,44.33){\line(1,5){7.00}} \put(2.67,43.67){\line(1,1){32.67}}

\end{picture}

\vskip -1cm

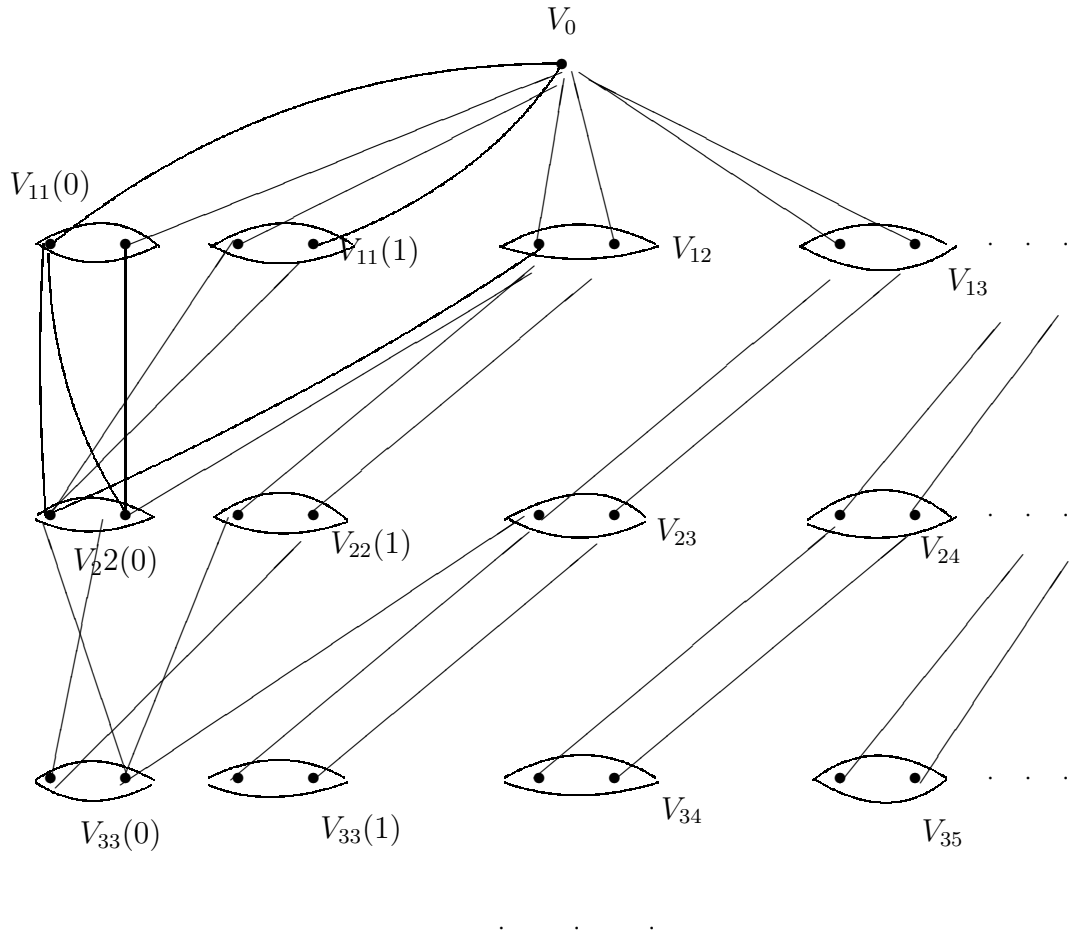
\centerline{Figure 1.} \vskip 1cm

In the proof of Theorem \ref{T4.2} we used the fact that a given Borel
dynamical system $(X,T)$ together with a vanishing sequence of markers can
be topologized by a topology $\om$ (without changing the Borel structure)
such that $X$ becomes a 0-dimensional Polish space, all atoms of partitions
defined by the sequence of markers are clopen, and $T$ becomes a
homeomorphism. The next theorem shows that if additionally the space $X$ is
{\it locally compact} in $\omega$, then applying the construction given in
subsection {\bf 4.2} we get a special Borel-Bratteli diagram.

\begin{theorem}\label{T4.4} Let $X$ and $T$ satisfy the made above assumptions.
Then $T$ is homeomorphic to a Vershik transformation acting on the space of
infinite paths of a special Borel-Bratteli diagram.
\end{theorem}

\noindent {\it Proof}. Let $(A_n)$ be a vanishing sequence of markers for
$T$ and let $\xi_n,\ n\in \N,$ be a partition of $X$ into disjoint
$T$-towers $\xi_n(v)= (T^j(A_n(v) : j=0,...,h(n,v)-1), v\in V_n$
constructed as in {\bf 4.2}. Our assumptions imply that $X$ can be taken to
be a locally compact 0-dimensional Polish space, $T$ is a homeomorphism,
and each atom $T^j(A_n(v))$ of $\xi_n,\ n\in \N$ is a clopen set.
\\

\noindent{\bf Claim 1} Let $X=\bigcup_{i\geq 0} X_i$ be a partition into
compact clopen disjoint sets $X_i$, then $B_n := \bigcup_{i\geq n} X_i,\
n=0,1,...,$ is a vanishing sequence of markers for $T$.
\\

\noindent{\it Proof}. To see this, we note that $A_n,\ n\in N$, is clopen
because the complement
$$
 A_n^c = X\setminus A_n = \bigcup_{v\in V_n}\bigcup_{j=1}^{h(n,v)-1}
T^j(A_n(v))
$$
is open. Then for any compact set $Y$, we have that $|\{n : A_n\cap Y\neq
\emptyset\}| <\infty$, that is $A_n \subset Y^c$ for sufficiently large $n$
(recall that the $A_n \supset A_{n+1}$ and $\bigcap_n A_n=\emptyset$). Given
$i$, choose $n_i\in \N$ such that $(\bigcup_{j=1}^iX_j)\cap A_n =\emptyset$
for $n \geq n_i$. Then $A_{n_i}\subset X_{i+1} \cup X_{i+2}\cup\cdots$,
hence $B_n$ is a complete $T$-section. Thus, $(B_n)$ is a vanishing sequence
of markers. It is easy to see that the partition $\xi_n$ (we do not change
our notation here), defined by $B_n$ and the homeomorphism $T$, also
consists of clopen sets for all $n$. The claim is proved.
\\

\noindent{\bf Claim 2} There exists a decomposition of $X$ into clopen
compact sets $X_0,X_1,X_2,...$ such that partitions $\xi_n,\ n\in \N,$
constructed by the vanishing sequence of markers $A_n := X_n\cup X_{n+1}
\cup\cdots,\ n=0,1,...$, have the following properties:
\begin{equation}\label{4.5}
h(n,v)=1\ \ {\rm if}\ A_n(v) \subset \bigcup_{i> n}X_i,
\end{equation}
\begin{equation}\label{4.6}
X_n = \bigcup_{v\in V_{nn}}A_n(v)\ \ {\rm where}\ V_{nn}\subset V_n,\
|V_{nn}|< \infty,
\end{equation}
\begin{equation}\label{4.7}
\bigcup_{i=1}^{n-1}X_i =\bigcup_{v\in V_{nn}(0)}\bigcup_{j\geq 1}^{h(n,v)-1}
T^j(A_n(v))
\end{equation}
where $V_{nn}(0) = \{v\in V_{nn} : h(n,v) > 1\}$.
\\

\noindent{\it Proof}. Take a partition of $X$ into compact clopen sets
$(Y_0,Y_1,Y_2,...)$. Denote $B_n := \bigcup_{i\geq n} Y_i,\ n\in \N,$ and
let $T^j(B_n(v)),\ j=0,1,...,h(n,v)-1, v\in V_n$, be elements of the
partition $\eta_n$ of $X$ into clopen $T$-towers constructed by $T$ and
$B_n$. Without loss of generality, we can assume that the $T$-towers
$(T^j(B_n(v)): j=0,1,...,h(n,v)-1)$ are chosen such that
\begin{equation}\label{4.8}
\lim_{n\to\infty}[\sup_{v\in V_n}( {\rm diam}\ T^j(B_n(v)))] =0.
\end{equation}

Consider these $T$-towers over $B_n,\ n\geq 1$. We can also assume that
every $Y_i$ is an $\eta_n$-set. Since $Y_0\cup Y_1\cup\cdots\cup Y_{n-1}$ is
compact, we can chose minimal $m=m(n)> n$ such that $h(n,v) = 1$ if
$B_n(v)\subset Y_i,\ i > m(n)$. Let $V_{nn}$ be a finite subset of $V_n$
such that the disjoint $T$-towers $\eta_n(v),\ v\in V_{nn}$, cover $Y_0\cup
Y_1\cup \cdots\cup Y_{m(n)}$. Some of those towers may be of height 1. Note
that the set $Y_0\cup\cdots\cup Y_{n-1}$ is covered by a finite subset of
$T$-towers $\eta_n(v)$ of heights greater than 1. Denote these towers by
$V_{nn}(0)$, that is
$$
Y_0\cup\cdots\cup Y_{n-1} = \bigcup_{v\in V_{nn}(0)}\bigcup_{j=
1}^{h(n,v)-1} T^j(B_n(v)).
$$
Let $V_{nn}(1) = V_{nn} \setminus V_{nn}(0)$. Then $h(v,n)=1$ when $v\in
V_{nn}(1)$.

Define a sequence of clopen compact disjoint sets $(X_n)$ by induction. let
$X_0 = Y_0$ and set
$$
\begin{array}{lll}
X_1 &= &Y_1\cup\cdots\cup Y_{m_1},\ \ \ \ \ \ \ m_1 = m(1)\\
X_2&=&Y_{m_1+1}\cup\cdots\cup Y_{m_2},\ \ \ m_2= m(m_1)\\
\end{array}
$$
$$
\centerline{$\cdots\cdots\cdots\cdots\cdots\cdots\cdots\cdots$}
$$
$$
\ X_{k+1}=\ Y_{m_k+1}\cup\cdots\cup Y_{m_{k+1}},\ \ m_{k+1}=m(m_k)
$$
$$
\centerline{$\cdots\cdots\cdots\cdots\cdots\cdots\cdots\cdots$}
$$

Let $A_n := X_n\cup X_{n+1}\cup\cdots,\ n\in \N$. Using $(A_n)$ and $T$,
construct refining partitions $(\xi_n)$, satisfying (\ref{4.8}). It follows
from the construction that $(\xi_n)$ is a telescoping of $(\eta_n)$. Then
$X$ is partitioned by $\xi_n$ into disjoint $T$-towers $\xi_n(v) =
(T^j(A_n(v)) : j=0,1,...,h(n,v)-1),\ v\in V_n$, such that conditions
(\ref{4.5}) - (\ref{4.7}) hold. The decomposition of $X= X_0\cup
X_1\cup\cdots$ and towers $\xi_n(v)$ are shown in Figure 2. The claim is
proved.
\\
\vskip 1cm

\unitlength=1.00mm \special{em:linewidth 0.4pt} \linethickness{0.4pt}
\begin{picture}(140.00,100.00)(5,45)
\linethickness{0.4pt} \put(130.33,80.00){\line(0,0){0.00}}
\linethickness{1.4pt} \put(9.67,80.33){\line(1,0){120.67}}
\linethickness{0.4pt} \put(19.67,81.67){\line(0,-1){2.33}}
\put(59.67,82.33){\line(0,-1){5.00}} \put(96.67,82.33){\line(0,-1){4.33}}
\put(132.67,80.00){\circle*{1.00}} \put(137.33,80.00){\circle*{1.00}}
\put(141.33,80.00){\circle*{1.00}} \put(10.33,82.00){\line(0,-1){3.67}}
\put(10.33,78.33){\line(0,0){0.00}} \put(10.33,82.67){\line(0,-1){5.00}}
\put(10.33,82.00){\line(0,-1){4.00}} \put(29.00,81.67){\line(0,-1){2.67}}
\put(39.33,81.67){\line(0,-1){2.67}} \linethickness{1.4pt}
\put(10.00,90.00){\line(1,0){9.67}} \put(28.67,90.00){\line(1,0){10.67}}
\put(49.33,90.00){\line(1,0){10.33}} \put(10.33,99.67){\line(1,0){8.67}}
\linethickness{1.4pt} \put(28.67,100.00){\line(1,0){10.00}}
\put(49.67,100.00){\line(1,0){10.00}} \linethickness{0.4pt}
\put(14.33,107.00){\circle*{1.00}} \put(33.67,107.00){\circle*{1.00}}
\put(55.00,107.00){\circle*{1.00}} \put(55.00,106.67){\circle*{1.00}}
\put(14.33,114.00){\circle*{1.00}} \put(33.67,114.33){\circle*{1.00}}
\put(55.33,115.00){\circle*{1.00}} \put(14.33,121.67){\circle*{1.00}}
\put(33.67,122.00){\circle*{1.00}} \put(56.00,122.33){\circle*{1.00}}
\linethickness{1.4pt} \put(10.67,129.67){\line(1,0){9.00}}
\put(28.67,129.67){\line(1,0){10.00}} \put(49.33,129.67){\line(1,0){10.33}}
\linethickness{0.4pt} \bezier{208}(59.67,134.00)(68.67,110.33)(58.33,86.00)
\bezier{204}(10.67,134.00)(1.00,110.00)(9.67,86.00)
\bezier{196}(11.00,134.67)(36.33,138.67)(59.00,134.33)
\bezier{196}(9.33,85.67)(32.33,85.00)(58.00,86.00)
\put(57.67,138.00){\makebox(0,0)[cc]{$X_0\cup\cdots\cup X_{n-1}$}}
\put(70.33,81.67){\line(0,-1){2.33}} \put(82.67,81.67){\line(0,-1){2.33}}
\put(107.67,81.67){\line(0,-1){2.33}} \put(119.00,81.33){\line(0,-1){2.00}}
\put(7.67,73.00){\makebox(0,0)[cc]{$X_n$}}
\put(64.00,73.00){\makebox(0,0)[cc]{$X_{n+1}$}}
\put(102.33,74.00){\makebox(0,0)[cc]{$X_{n+2}$}}
\put(18.67,60.33){\makebox(0,0)[cc]{$V_{nn}(0)$}}
\put(43.33,59.33){\makebox(0,0)[cc]{$V_{nn}(1)$}}
\put(42.33,63.00){\vector(-1,1){16.33}}
\put(42.00,63.33){\vector(1,4){3.67}} \put(18.00,63.67){\vector(-1,4){3.67}}
\put(18.33,63.00){\vector(1,1){16.00}}
\put(18.33,63.67){\vector(3,1){35.67}}
\put(74.67,85.67){\makebox(0,0)[cc]{$A_n(v)$}}
\put(49.33,82.00){\line(0,-1){2.67}} \put(49.33,81.67){\line(0,-1){3.00}}
\end{picture}
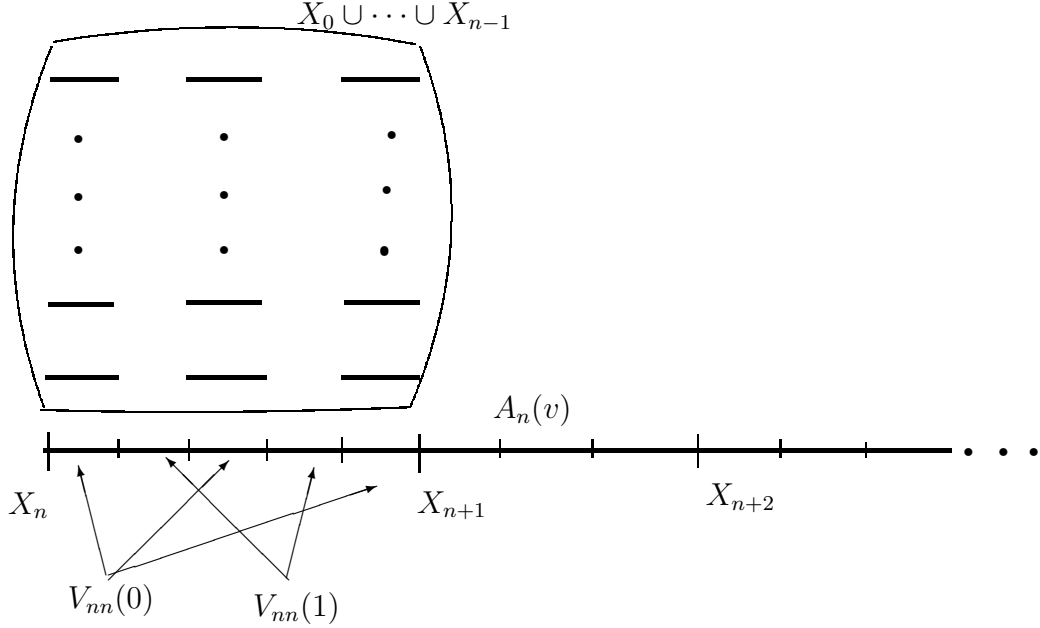
\vskip -0.5cm \centerline{Figure 2} \vskip 1cm

To complete the proof of the theorem, we need to show that the ordered
Borel-Bratteli diagram $B=(V,E,\geq)$ constructed by the vanishing sequence
of markers $(A_n)$ defined in the proof of Claim 2 satisfies properties (i)
and (ii) of the definition of a special Borel-Bratteli diagram.

Denote by $V_{nj} = \{v\in V_n : A_n(v) \subset X_j\},\ j\geq n, n\in \N$.
Let $E_{nj}$ be the set of edges between $V_{nj}$ and $V_{n-1,j}$. Then
$V_n$ is partitioned into non-empty finite sets $V_{nj},\ j\geq n$. Clearly,
we can assume that $|V_{nj}| \geq 2$ for all $n,j$ refining, in case of
need, the partition $\xi_n$. It is obvious that for every $A_n(v),\ v\in
V_{nj},\ j > n,$ there exists a unique set $A_{n-1}(v'),\ v'\in V_{n-1,j}$
such that $A_n(v) \subset A_{n-1}(v')$. This means that the set $E(v',v)$ of
edges connecting $v'$ and $v$ has exactly one element.

The set $V_{nn}$ is divided into two sets $V_{nn}(0) = \{v\in V_{nn} :
h(n,v) > 1\}$ and $V_{nn}(1) = \{v\in V_{nn} : h(n,v) =1\},\ n\in \N$. Then
\begin{equation}\label{4.9}
X_0\cup\cdots\cup X_{n-1} = \bigcup_{v\in V_{nn}(0)}
\bigcup_{j=1}^{h(n,v)-1} T^j(A_n(v))
\end{equation}
\begin{equation}\label{4.10}
X_n\cap T^{-1}(X_n\cup X_{n+1}\cup\cdots) = \bigcup_{v\in V_{nn}(1)} A_n(v).
\end{equation}
Moreover,
\begin{equation}\label{4.11}
 \bigcup_{v\in V_{nn}(0)} A_n(v) = X_n\cap T^{-1}(X_0\cup X_{1}\cup\cdots
\cup X_{n-1}) = X_n\cap T^{-1}(X_{n-1}).
\end{equation}
Relations (\ref{4.9}) - (\ref{4.11}) show that $r^{-1}(v)$ consists of a
single edge connecting $v\in V_{nn}(1)$ and a vertex $w\in V_{n-1,n}$.
Without loss of generality, we can assume that $|V_{nn}(0)| \geq 2$. If
$v\in V_{nn}(0)$, then the set $r^{-1}(v)$ has a unique edge $e(w,v)$
connecting $v$ with some $w\in V_{n-1,n}$ and a number of edges which
connect $v$ with vertices from $V_{n-1,n-1}$. It follows from the
construction of refining sequence $(\xi_n)$ that $s(E_{nj}) = V_{n-1,j}$ for
$j>n$. We get also from (\ref{4.11}) that $s(E_{nn}(0)) = V_{n-1,n-1}$ where
$E_{nn}(0)$ is the set of edges arriving to vertices from $V_{nn}(0)$. A
non-trivial linear order on $r^{-1}(v)$ should be defined for $v\in
V_{nn}(0)$. We assign the least value 1 to the edge $e(w,v)$. The maximal
value of the defined order on $r^{-1}(v)$ is assigned to an edge $e$ such
that $s(e)$ is either in $V_{n-1,n}$ or in $V_{n-1,n-1}(1)$. It is obvious
that the set $Y_B$ of infinite paths is uncountable and has no cofinal
maximal and minimal points. The fact that $T$ is isomorphic to the Vershik
automorphism acting on $Y_B$ is proved as in Theorem \ref{T4.2} and
Proposition \ref{P4.4}. \hfill$\square$
\\

{\it Acknowledgement}. We are grateful to many people  for discussing our
results. We would like to express our special thanks to A.~Kechris,
G.~Hjorth, and B.~Miller for numerous talks which helped us to improve the
first version of the paper. The first-named author thanks the University of
New South Wales for the warm hospitality and the Australian Research
Council for its support. S.B. and J.K. acknowledge also the support of the
Torun University and the Kharkov Institute for Low Temperature Physics
during exchange visits.


\vskip1cm \noindent {\small {\it S.~Bezuglyi\\
Department of Mathematics\\
Institute for Low Temperature Physics\\
47, Lenin ave.\\
61103 Kharkov, UKRAINE\\
\smallskip
\noindent bezuglyi@ilt.kharkov.ua}}
\\
\\
\noindent {\small {\it A.H.~Dooley\\
School of Mathematics\\
University of New South Wales\\
Sydney, 2052 NSW, AUSTRALIA\\
\smallskip
\noindent tony@maths.unsw.edu.au}}
\\
\\
\noindent {\small {\it J.~Kwiatkowski\\
The College of Economics and Computer Sciences\\
Wyzwolenia 30\\
10 - 106, Olsztyn, POLAND\\
\smallskip
\noindent jkwiat@mat.uni.torun.pl}}

\end{document}